# Optimal Pairs Trading with Time-Varying Volatility

T. N. Lǐ* and A. Tourin†

January 6, 2016

**Abstract**

We propose a pairs trading model that incorporates a time-varying volatility of the Constant Elasticity of Variance type. Our approach is based on stochastic control techniques; given a fixed time horizon and a portfolio of two co-integrated assets, we define the trading strategies as the portfolio weights maximizing the expected power utility from terminal wealth. We compute the optimal pairs strategies by using a Finite Difference method. Finally, we illustrate our results by conducting tests on historical market data at daily frequency. The parameters are estimated by the Generalized Method of Moments.

## 1 Introduction

This article extends the pairs trading model proposed in Tourin and Yan (2013) to incorporate time-varying volatility. In contrast to Tourin and Yan (2013), we are unable to derive a fully explicit solution for the optimal pairs trading strategies, even for a fully specified local volatility model. However, by using stochastic control techniques, we show that this problem can be reduced to a degenerate two-dimensional linear Partial Differential Equation (PDE in short) and can therefore easily be solved numerically by using an implicit Finite Difference scheme. We refer to the book by Fleming and Soner (1993) and Pham (2009) for an introduction to stochastic control and its applications in Finance, and to Touzi (2013) for a presentation of Finite Difference schemes for financial applications.

The volatility model we choose is of the Constant Elasticity of Variance type, introduced in Cox (1975), Cox and Ross (1976), which captures the negative link between volatility and stock price observed on financial data, for a relatively modest computational cost.

Although the solution can be easily computed for both the power and the exponential utility function, we mainly adopt the power utility function in this paper because it allows us to compute a rate of return for the computed strategies. However, for a sake of completeness, we briefly present in section 2 the results for the exponential utility function in the constant volatility case.

Furthermore, we illustrate our approach with historical market data at daily frequency. We estimate the parameters by using the Generalized Method of Moments introduced in Hansen (1982),

*Department of Mathematics, New York University, 251 Mercer Street, New York, New York 10012. *nl747@cims.nyu.edu*

†Department of Finance and Risk Engineering, New York University, Six Metrotech Center, Brooklyn, NY 11201 *atourin@nyu.edu.*

and compute numerically the associated optimal trading strategies by using the Finite Difference scheme.

A body of literature on the application of stochastic control theory to pairs trading has emerged in the past several years, starting with the work of Mudchanatongsuk et al. (2008). Recently, Leung and Li (2015) formulated a model characterizing the optimal entry and exit points of a pairs trading strategy under transaction costs. Then, Lei and Xu (2015) proposed a model for determining multiple entry and exit-points during a trading period, and illustrated their results with applications to dual-listed Chinese stocks. Ngo and Pham (2014) frames the pairs trading problem as a regime switching model between three regimes: flat positions, one long position on one asset and a short position on the other, and vice versa . Finally, Lintilhac (2014), Cartea and Jaimungal (2015) considered applications to portfolios of co-integrated assets, both generalizing the model in Tourin and Yan (2013). Lintilhac (2014) illustrated the model by applying it to historical data in the bitcoin markets while Cartea and Jaimungal (2015) performed simulations to analyse its performance.

This paper is organized as follows. First, we revisit in the next section the constant volatility case and solve it explicitly in a more general setting and for a different utility function than in Tourin and Yan (2013). Specifically, we consider a non-zero risk-free rate, add a deterministic trend in the co-integration variable, make the dynamics fully symmetric by incorporating the co-integration vector in the drifts of both stocks, and use the power utility function. We also provide a new verification Theorem. The main purpose of this first section is to introduce through the description of the pairs trading problem in the constant volatility case, the building blocks for the more general model with stochastic volatility. It also serves as a benchmark the numerical Finite Difference scheme can be tested on.

In the third section, we present the main model with a time-varying local volatility, derive the Hamilton-Jacobi-Bellman equation characterizing the optimal trading strategies by using standard stochastic control techniques, and reduce it to a two-dimensional linear degenerate parabolic PDE, showing that the control can be decomposed into the sum of an time-independent component, which can be written explicitly in terms of the parameters in the model, and a time-varying component which must be computed approximately when the volatility is no longer constant.

Finally, in section 4, we describe briefly the Finite Difference method we design to compute the pairs trading strategies and present some experiments for daily-frequency market data.

## 2 The model with constant volatility

In this section, we derive an exact solution in closed form for the pairs trading problem with constant volatility, and prove a verification Theorem ensuring that the computed solution is the solution of the stochastic optimal control problem. We choose the power utility function despite the fact that it makes the computations more complex than for the exponential utility function because, as we explain below, the power utility function possesses a major advantage: it allows us to compute a realized rate of return when backtesting the strategies.

We essentially use for the co-integrated stocks the model derived in Duan and Pliska (2004) as the diffusion limit of a discrete-time model, in the case when there are only two assets and the volatility is constant. We fix a time horizon $T > 0$. The co-integrated asset prices $S_t^1$ and $S_t^2$, for

$t \in [0, T]$, satisfy the Stochastic Differential Equations

$$d \log S_t^1 = \left( \mu_1 - \frac{\sigma_1^2}{2} + \delta_1 z_t \right) dt + \sigma_1 dB_t^1 \tag{2.1}$$

$$d \log S_t^2 = \left( \mu_2 - \frac{\sigma_2^2}{2} + \delta_2 z_t \right) dt + \sigma_2 dB_t^2 \tag{2.2}$$

where $\mu_1, \mu_2, \delta_1, \delta_2, \sigma_1 > 0, \sigma_2 > 0$ are constant, $\left(B_t^1, B_t^2\right)$ is a pair of standard Brownian motions with correlation coefficient $\rho \in (-1, 1)$, and the co-integrating vector $z_t$ is defined by

$$z_t = a + bt + \log S_t^1 + \beta \log S_t^2. \tag{2.3}$$

We denote the underlying filtration generated by the pair of Brownian motions by $\mathcal{F}_t = \sigma\big( \left(B_s^1, B_s^2\right) : 0 \leq s \leq t\big)$. Furthermore, $z_t$ satisfies the dynamics

$$\begin{aligned} dz_t &= \left( b + \mu_1 - \frac{\sigma_1^2}{2} + \beta\mu_2 - \beta\frac{\sigma_2^2}{2} + \delta_1 z_t + \beta\delta_2 z_t \right) dt + \sigma_1 dB_t^1 + \beta\sigma_2 dB_t^2 \\ &= \alpha \left( \eta - z_t \right) dt + \sigma_\beta dB_t \end{aligned}$$

where $\alpha = -\delta_1 - \beta\delta_2$ is the speed of mean reversion and must satisfy $0 < \alpha < 2$ for the log-prices to be co-integrated (see Proposition 1 in Duan and Pliska (2004)), $\sigma_\beta = \sqrt{\sigma_1^2 + \beta^2\sigma_2^2}$, $B_t = \frac{\sigma_1}{\sigma_\beta} B_t^1 + \beta\frac{\sigma_2}{\sigma_\beta} B_t^2$ is a Brownian motion adapted to $\mathcal{F}_t$ and $\eta = \frac{1}{\alpha}\left( b + \mu_1 - \frac{\sigma_1^2}{2} + \beta\left( \mu_2 - \frac{\sigma_2^2}{2} \right) \right)$ is the equilibrium level.

We also assume that there is a risk free asset such as a money market account with interest rate $r > 0$. Its price per share is denoted by $S_t^0$, for $t \in [0, T]$ and evolves according to $dS_t^0 = rS_t^0 dt$.

Next, the variable $W_s$ represents the value of the investor's portfolio at time $s$. Starting at time $t$ with an initial wealth $W_t = w$, she invests, at every time $s \in [t, T]$, in both the risk-free money market account and the two stocks. We denote by $\pi_s^1, \pi_s^2$ the fractions of wealth invested respectively in the first and second stocks at time $s$. We only consider self-financing strategies and hence, the evolution of the wealth variable is given by

$$dW_s = \pi_s^1 W_s \frac{dS_s^1}{S_s^1} + \pi_s^2 W_s \frac{dS_s^2}{S_s^2} + \left(1 - \pi_s^1 - \pi_s^2\right) W_s \frac{dS_s^0}{S_s^0}.$$

We also consider the log-prices $x_s = \log S_s^1, y_s = \log S_s^2$. Finally, the dynamics of the state variables $W_s, x_s, y_s$ are given by the following controlled system of Stochastic Differential Equations

$$\begin{aligned} dW_s &= \pi_s^1 \left(\mu_1 + \delta_1 z_s\right) W_s ds + \pi_s^2 \left(\mu_2 + \delta_2 z_s\right) W_s ds \\ &\qquad + r(1 - \pi_s^1 - \pi_s^2) W_s ds + \pi_s^1 \sigma_1 W_s dB_s^1 + \pi_s^2 \sigma_2 W_s dB_s^2, \end{aligned} \tag{2.4}$$

$$dx_s = \left( \mu_1 - \frac{\sigma_1^2}{2} + \delta_1 z_s \right) ds + \sigma_1 dB_s^1 \tag{2.5}$$

$$dy_s = \left( \mu_2 - \frac{\sigma_2^2}{2} + \delta_2 z_s \right) ds + \sigma_2 dB_s^2, \tag{2.6}$$

$$W_t = w,\ x_t = x,\ y_t = y, \tag{2.7}$$

where $z_t$ is defined in (2.3). Furthermore, we require the wealth to be non-negative, i.e.

$$w_s \geq 0 \text{ a.s. for all } s \in [0, T]\,. \tag{2.8}$$

A pair of controls $\left(\pi^1, \pi^2\right)$ is said to be admissible if $\pi^1, \pi^2$ are real-valued, progressively measurable, are such that, (2.4), (2.5), (2.6), and (2.7) define a unique solution $(W_s, x_s, y_s)$ for every time $s \in [0, T]$, the wealth variable satisfies the constraint (2.8), and $\left(\pi^1, \pi^2, W_s\right)$ satisfy the integrability condition

$$\mathbb{E}\int_t^T \left(\pi_s^1 W_s\right)^2 + \left(\pi_s^2 W_s\right)^2 ds < +\infty.$$

We denote the set of admissible controls at the initial time of investment $t$, by $\mathcal{A}_t$. Next, we define the value function $u\left(t, w, x, y\right)$ of the following backward dynamic optimization problem: the investor seeks an admissible strategy $\left(\pi_s^1, \pi_s^2\right)$ for every $s \in [t, T)$, that maximizes the utility he derives from wealth at time $T$, i.e.

$$u\left(t, w, x, y\right) = \sup_{\left(\pi^1, \pi^2\right) \in \mathcal{A}_t} \mathbb{E}\left[U\left(W_T^{t,w,x,y,\left(\pi^1,\pi^2\right)}\right)\right], \tag{2.9}$$

where $W_T^{t,w,x,y,\left(\pi^1,\pi^2\right)}$ denotes the solution of (2.4) at time $T$, corresponding to the control pair $\left(\pi^1, \pi^2\right)$ and the initial conditions (2.7).

Furthermore, we focus on the case of the power utility function, i.e.

$$U\left(w\right) = \frac{1}{\gamma} w^\gamma, \tag{2.10}$$

where $\gamma > 0$ denotes the constant risk aversion coefficient.

In order to solve this problem, we apply standard stochastic control techniques (see for instance Fleming and Soner (1993)). Since we do not know a priori the regularity of the value function of the stochastic control problem, we only proceed formally with the hope of obtaining a smooth candidate solution for the stochastic control problem. We will verify later that our candidate solution coincides with the solution of the optimal stochastic control. We expect the value function $u$ defined in (2.9) to satisfy the following Hamilton-Jacobi-Bellman (HJB in short) equation:

$$\begin{aligned} &- u_t - \sup_{\pi_1, \pi_2} \Bigg\{ \left[\pi_1\left(\mu_1 + \delta_1 z\right) + \pi_2\left(\mu_2 + \delta_2 z\right)\right] w u_w + r\left(1 - \pi_1 - \pi_2\right) w u_w \\ &+ \left(\mu_1 - \frac{1}{2}\sigma_1^2 + \delta_1 z\right) u_x + \left(\mu_2 - \frac{1}{2}\sigma_2^2 + \delta_2 z\right) u_y \\ &+ \pi_1 \sigma_1^2 w u_{wx} + \pi_2 \rho \sigma_1 \sigma_2 w u_{wx} + \pi_2 \sigma_2^2 w u_{wy} + \pi_1 \rho \sigma_1 \sigma_2 u_{wy} \\ &+ \frac{1}{2}\left(\pi_1^2 \sigma_1^2 + \pi_2^2 \sigma_2^2 + 2\pi_1 \pi_2 \rho \sigma_1 \sigma_2\right) w^2 u_{ww} + \frac{1}{2}\sigma_1^2 u_{xx} + \frac{1}{2}\sigma_2^2 u_{yy} + \rho \sigma_1 \sigma_2 u_{xy} \Bigg\} = 0, \end{aligned} \tag{2.11}$$

for all $0 \leq t < T$, $w \geq 0$, $x$, and $y$, coupled with the final condition

$$\begin{aligned} u\left(T, w, x, y\right) &= U\left(w\right) \\ &= \frac{1}{\gamma} w^\gamma, \quad \text{for all } w \geq 0,\ x,\ y. \end{aligned}$$

Next, with the goal of factoring out the wealth variable we define the ansatz

$$u\left(t, w, x, y\right) = \frac{1}{\gamma} w^\gamma g\left(t, x, y\right),$$

and we derive the PDE for $g(t,x,y)$

$$\begin{aligned}
&-g_t-\sup_{\pi_1,\pi_2}\Bigg\{[\pi_1(\mu_1+\delta_1 z)+\pi_2(\mu_2+\delta_2 z)]\gamma g+r(1-\pi_1-\pi_2)\gamma g\\
&+\left(\mu_1-\frac{1}{2}\sigma_1^2+\delta_1 z\right)g_x+\left(\mu_2-\frac{1}{2}\sigma_2^2+\delta_2 z\right)g_y\\
&+\left(\pi_1\sigma_1^2+\pi_2\rho\sigma_1\sigma_2\right)\gamma g_x+\left(\pi_2\sigma_2^2+\pi_1\rho\sigma_1\sigma_2\right)\gamma g_y\\
&+\frac{1}{2}\left(\pi_1^2\sigma_1^2+\pi_2^2\sigma_2^2+2\pi_1\pi_2\rho\sigma_1\sigma_2\right)\gamma(\gamma-1)g+\frac{1}{2}\sigma_1^2 g_{xx}+\frac{1}{2}\sigma_2^2 g_{yy}+\rho\sigma_1\sigma_2 g_{xy}\Bigg\}=0,
\end{aligned}$$

for all $0\le t<T$, $x$, and $y$, coupled with $g(T,x,y)=1$ for all $x$ and $y$.

Finally, we can replace the pair $(x,y)$ by the single variable $z$ in order to reduce the dimension further. Unfortunately this will not work when the volatility coefficients are not constant as we will see in the next section. To this end, we define $h(t,z)=g(t,x,y)$, and derive the one-dimensional HJB equation satisfied by $h(t,z)$

$$\begin{aligned}
&-h_t-\sup_{\pi_1,\pi_2}\Bigg\{[\pi_1(\mu_1+\delta_1 z)+\pi_2(\mu_2+\delta_2 z)]\gamma h+r(1-\pi_1-\pi_2)\gamma h\\
&+\left[b+\left(\mu_1-\frac{1}{2}\sigma_1^2+\delta_1 z\right)+\beta\left(\mu_2-\frac{1}{2}\sigma_2^2+\delta_2 z\right)\right]h_z+\big[\pi_1\sigma_1^2\\
&+\pi_2\sigma_2^2\beta+\pi_2\rho\sigma_1\sigma_2+\pi_1\beta\rho\sigma_1\sigma_2\big]\gamma h_z\\
&+\frac{1}{2}\left(\pi_1^2\sigma_1^2+\pi_2^2\sigma_2^2+2\pi_1\pi_2\rho\sigma_1\sigma_2\right)\gamma(\gamma-1)h+\frac{1}{2}\left(\sigma_1^2+\sigma_2^2\beta^2+2\beta\rho\sigma_1\sigma_2\right)h_{zz}\Bigg\}=0,
\end{aligned}$$

for all $0\le t<T$ and $z$, coupled with $h(T,z)=1$ for all $z$.

Next, the optimal controls are given in terms of $h$ and its partial derivatives by

$$\pi_1^*=\frac{\mu_1-r+\delta_1 z}{\sigma_1^2(1-\gamma)(1-\rho^2)}-\rho\frac{\mu_2-r+\delta_2 z}{\sigma_1\sigma_2(1-\gamma)(1-\rho^2)}+\frac{h_z}{(1-\gamma)h} \tag{2.12}$$

$$\pi_2^*=\frac{\mu_2-r+\delta_2 z}{\sigma_2^2(1-\gamma)(1-\rho^2)}-\rho\frac{\mu_1-r+\delta_1 z}{\sigma_1\sigma_2(1-\gamma)(1-\rho^2)}+\frac{\beta h_z}{(1-\gamma)h}. \tag{2.13}$$

Substituting the controls (2.12) and (2.13) into the PDE yields

$$\begin{aligned}
&-h_t+\frac{\gamma}{2(\gamma-1)(1-\rho^2)} \qquad (2.14)\\
&\left[\frac{(\mu_1-r+\delta_1 z)^2}{\sigma_1^2}+\frac{(\mu_2-r+\delta_2 z)^2}{\sigma_2^2}-2\frac{\rho(\mu_1-r+\delta_1 z)(\mu_2-r+\delta_2 z)}{\sigma_1\sigma_2}\right]h\\
&-r\gamma h+\frac{1}{\gamma-1}[(\mu_1+\delta_1 z)+\beta(\mu_2+\delta_2 z)]h_z\\
&-\frac{r\gamma(1+\beta)}{\gamma-1}h_z-bh_z+\frac{1}{2}\left(\sigma_1^2+\beta\sigma_2^2\right)h_z\\
&+\frac{\gamma}{2(\gamma-1)}\left(\sigma_1^2+\beta^2\sigma_2^2+2\beta\rho\sigma_1\sigma_2\right)\frac{h_z^2}{h}-\frac{1}{2}\left(\sigma_1^2+\sigma_2^2\beta^2+2\beta\rho\sigma_1\sigma_2\right)h_{zz}=0,
\end{aligned}$$

with the terminal condition $h(T,z)=1$.

We can remove the non-linearity in (2.14) by making the change of unknown function

$$h=\frac{1}{1-\gamma}\phi^{1-\gamma}.$$

The linear PDE for $\phi$ is

$$
\begin{aligned}
&-\phi_t-\frac{\gamma}{2(\gamma-1)^2(1-\rho^2)} \\
&\left[\frac{(\mu_1+\delta_1 z-r)^2}{\sigma_1^2}+\frac{(\mu_2+\delta_2 z-r)^2}{\sigma_2^2}-2\frac{\rho(\mu_1-r+\delta_1 z)(\mu_2-r+\delta_2 z)}{\sigma_1\sigma_2}\right]\phi \\
&-\frac{r\gamma}{1-\gamma}\phi+\frac{1}{\gamma-1}\left[(\mu_1+\delta_1 z)+\beta(\mu_2+\delta_2 z)\right]\phi_z-\frac{r\gamma(1+\beta)}{\gamma-1}\phi_z \\
&-b\phi_z+\frac{1}{2}\left(\sigma_1^2+\beta\sigma_2^2\right)\phi_z-\frac{1}{2}\left(\sigma_1^2+\beta^2\sigma_2^2+2\beta\rho\sigma_1\sigma_2\right)\phi_{zz}=0,
\end{aligned}
\tag{2.15}
$$

with the terminal condition $\phi(T,z)=(1-\gamma)^{\frac{1}{1-\gamma}}$.

Finally, as in Benth and Karlsen (2005), we solve (2.15) by using the ansatz function $\phi(t,z)=\exp\{f_2(t)z^2+f_1(t)z+f_0(t)\}$ and deriving the Ordinary Differential Equations (ODEs in short) satisfied by the coefficients $f_2(t)$, $f_1(t)$, $f_0(t)$. The fist ODE, obtained from the coefficient of the term $z^2$ in (2.15) can be written as

$$
\begin{aligned}
&f_2'(t)+2\left(\sigma_1^2+\beta^2\sigma_2^2+2\beta\rho\sigma_1\sigma_2\right)\left(f_2(t)-\frac{\delta_1+\beta\delta_2}{2(\sigma_1^2+\beta^2\sigma_2^2+2\beta\rho\sigma_1\sigma_2)(\gamma-1)}\right)^2 \\
&-\frac{(\delta_1+\beta\delta_2)^2}{2(\gamma-1)^2(\sigma_1^2+\beta^2\sigma_2^2+2\beta\rho\sigma_1\sigma_2)}+\frac{\gamma}{2(\gamma-1)^2(1-\rho^2)}\left(\frac{\delta_1^2}{\sigma_1^2}+\frac{\delta_2^2}{\sigma_2^2}-2\rho\frac{\delta_1\delta_2}{\sigma_1\sigma_2}\right)=0
\end{aligned}
\tag{2.16}
$$

coupled with the terminal condition $f_2(T)=0$. The second ODE, corresponding to the coefficient of the term $z$ in (2.15) is

$$
\begin{aligned}
&f_1'(t)+\left[-\frac{1}{\gamma-1}(\delta_1+\beta\delta_2)+2\left(\sigma_1^2+\beta^2\sigma_2^2+2\beta\rho\sigma_1\sigma_2\right)f_2(t)\right]f_1(t) \\
&+\left[2b+\frac{2r\gamma(1+\beta)}{\gamma-1}-\frac{1}{\gamma-1}[2(\mu_1+\beta\mu_2)]-\left(\sigma_1^2+\beta\sigma_2^2\right)\right]f_2(t) \\
&+\frac{\gamma}{2(\gamma-1)^2(1-\rho^2)}\left[\frac{(\mu_1-r)\delta_1}{\sigma_1^2}+\frac{(\mu_2-r)\delta_2}{\sigma_2^2}-2\rho\left(\frac{\delta_1(\mu_2-r)}{\sigma_1\sigma_2}+\frac{\delta_2(\mu_1-r)}{\sigma_1\sigma_2}\right)\right]=0,
\end{aligned}
\tag{2.17}
$$

and must be coupled with the terminal condition $f_1(T)=0$. The third ODE, derived from the constant coefficient in (2.15) is

$$
\begin{aligned}
&f_0'(t)+\frac{\gamma}{2(\gamma-1)^2(1-\rho^2)}\left[\frac{(\mu_1-r)^2}{\sigma_1^2}+\frac{(\mu_2-r)^2}{\sigma_2^2}-2\rho\frac{(\mu_1-r)(\mu_2-r)}{\sigma_1\sigma_2})\right] \\
&+\frac{r\gamma}{\gamma-1}-\frac{1}{\gamma-1}(\mu_1+\beta\mu_2)f_1(t)+bf_1(t)-\frac{1}{2}\left(\sigma_1^2+\beta\sigma_2^2\right)f_1(t)+r\gamma\frac{(1+\beta)}{\gamma-1}f_1(t) \\
&+\frac{1}{2}\left(\sigma_1^2+\beta^2\sigma_2^2+2\beta\rho\sigma_1\sigma_2\right)f_1^2(t)+\left(\sigma_1^2+\beta^2\sigma_2^2+2\beta\rho\sigma_1\sigma_2\right)f_2(t)=0,
\end{aligned}
\tag{2.18}
$$

with the terminal condition $f_0(T)=\frac{1}{1-\gamma}\ln(1-\gamma)$.

For convenience, we define the constants

$$
\begin{aligned}
c_1&=\sigma_1^2+\beta^2\sigma_2^2+2\beta\rho\sigma_1\sigma_2>0,\\
c_2&=\frac{\alpha}{2(1-\gamma)c_1}>0,
\end{aligned}
$$

$$c_0 = \frac{\alpha^2}{2(1-\gamma)^2 c_1} - \frac{\gamma}{2(1-\gamma)^2(1-\rho^2)}\left(\frac{\delta_1^2}{\sigma_1^2} + \frac{\delta_2^2}{\sigma_2^2} - 2\rho\frac{\delta_1\delta_2}{\sigma_1\sigma_2}\right).$$

The solution $f_2$ depends on the sign of $c_0$. We have the result

**Lemma 2.1.**
- *If $c_0 > 0$, then*

$$f_2(t) = c_2\left(1 - \frac{c_0}{2c_1c_2^2}\right)\frac{\sinh\left(\sqrt{2c_1c_0}(T-t)\right)}{\sinh\left(\sqrt{2c_1c_0}(T-t)\right) + \frac{1}{c_2}\sqrt{\frac{c_0}{2c_1}}\cosh\left(\sqrt{2c_1c_0}(T-t)\right)}.$$

- *If $c_0 < 0$, then*

$$f_2(t) = c_2\left(1 - \frac{c_0}{2c_1c_2^2}\right)\frac{\sin\left(\sqrt{-2c_1c_0}(T-t)\right)}{\sin\left(\sqrt{-2c_1c_0}(T-t)\right) + \frac{1}{c_2}\sqrt{\frac{-c_0}{2c_1}}\cos\left(\sqrt{-2c_1c_0}(T-t)\right)}.$$

- *If $c_0 = 0$, then*

$$f_2(t) = c_2 + \frac{c_2}{-2c_1c_2(T-t) - 1}.$$

*Proof.* It is easy to verify that the above solution satisfies (2.16). □

It is then straightforward to deduce the solutions $f_0(t), f_1(t)$ of the other two ODEs in terms of $f_2$:

**Lemma 2.2.** *1. The solution $f_0$ of (2.18) is given by*

$$\begin{aligned} f_0(t) = & \int_t^T \left[\frac{1}{1-\gamma}(\mu_1 + \beta\mu_2) f_1(s) + b f_1(s) - \frac{1}{2}\left(\sigma_1^2 + \beta\sigma_2^2\right) f_1(s) \right. \\ & - r\gamma\frac{(1+\beta)}{1-\gamma} f_1(s) + \frac{1}{2}\left(\sigma_1^2 + \beta^2\sigma_2^2 + 2\beta\rho\sigma_1\sigma_2\right) f_1^2(s) \\ & \left. + \left(\sigma_1^2 + \beta^2\sigma_2^2 + 2\beta\rho\sigma_1\sigma_2\right) f_2(s)\right] ds + \frac{\gamma}{2(1-\gamma)^2(1-\rho^2)} \\ & \left[\frac{(\mu_1 - r)^2}{\sigma_1^2} + \frac{(\mu_2 - r)^2}{\sigma_2^2} - 2\rho\frac{(\mu_1 - r)(\mu_2 - r)}{\sigma_1\sigma_2}\right](T-t) \\ & - \frac{r\gamma}{1-\gamma}(T-t) + \frac{1}{1-\gamma}\ln(1-\gamma). \end{aligned}$$

*2. The solution $f_1$ of (2.17) is given by*

$$\begin{aligned} f_1(t) &= \int_t^T b_1(s)\, exp\{\int_t^s a_1(u)\, du\} ds, \ \textit{where} \\ a_1(t) &= \frac{-\alpha}{1-\gamma} + 2c_1 f_2(t), \ \textit{and} \\ b_1(t) &= \frac{\gamma}{2(1-\gamma)^2(1-\rho^2)} \\ & \left[\frac{(\mu_1 - r)\delta_1}{\sigma_1^2} + \frac{(\mu_2 - r)\delta_2}{\sigma_2^2} - 2\rho\left(\frac{\delta_1(\mu_2 - r)}{\sigma_1\sigma_2} + \frac{\delta_2(\mu_1 - r)}{\sigma_1\sigma_2}\right)\right] \\ & + \left[2b - \frac{2r\gamma(1+\beta)}{1-\gamma} + \frac{1}{1-\gamma}[2(\mu_1 + \beta\mu_2)] - \left(\sigma_1^2 + \beta\sigma_2^2\right)\right] f_2(t). \end{aligned}$$

**Lemma 2.3.** *The optimal controls are of the form*

$$\pi_1^* = \frac{\mu_1 - r + \delta_1 z}{\sigma_1^2 (1-\gamma)(1-\rho^2)} - \rho \frac{\mu_2 - r + \delta_2 z}{\sigma_1 \sigma_2 (1-\gamma)(1-\rho^2)} + 2 f_2(t) z + f_1(t) \tag{2.19}$$

$$\pi_2^* = \frac{\mu_2 - r + \delta_2 z}{\sigma_2^2 (1-\gamma)(1-\rho^2)} - \rho \frac{\mu_1 - r + \delta_1 z}{\sigma_1 \sigma_2 (1-\gamma)(1-\rho^2)} + \beta(2 f_2(t) z + f_1(t)). \tag{2.20}$$

*where $f_1$, $f_2$, and $f_3$ are defined in lemmas 2.1 and 2.2.*

Note that they are three components in these controls: the first two are time-independent, whereas the third one is time-dependent and vanishes at time $t = T$. Indeed this model is not time-consistent and we partially unwind the positions progressively over time as we approach the time horizon.

Furthermore, note that the rate of return over a small time interval is independent of the initial capital and is given by

$$\frac{dW_t}{W_t} = rdt + \pi_{1,t}^* \left( \frac{dS_t^1}{S_t^1} - rdt \right) + \pi_{2,t}^* \left( \frac{dS_t^2}{S_t^2} - rdt \right), \tag{2.21}$$

where $\pi_{1,t}^*, \pi_{2,t}^*$ are the optimal control processes defined in Lemma 2.3, i.e.

$$\pi_{1,t}^* = \frac{\mu_1 - r + \delta_1 z_t}{\sigma_1^2 (1-\gamma)(1-\rho^2)} - \rho \frac{\mu_2 - r + \delta_2 z_t}{\sigma_1 \sigma_2 (1-\gamma)(1-\rho^2)} + 2 f_2(t) z_t + f_1(t)$$

$$\pi_{2,t}^* = \frac{\mu_2 - r + \delta_2 z_t}{\sigma_2^2 (1-\gamma)(1-\rho^2)} - \rho \frac{\mu_1 - r + \delta_1 z_t}{\sigma_1 \sigma_2 (1-\gamma)(1-\rho^2)} + \beta \left( 2 f_2(t) z_t + f_1(t) \right).$$

Next, we prove the following verification Theorem, which is very similar to the main Theorem in Benth and Karlsen (2005).

**Theorem 2.1.** *If the following conditions are satisfied*

1. *If*

$$(1-\gamma) f_2(t) (e^{2\alpha(T-t)} - 1) < \frac{\alpha}{2\sigma_\beta^2},$$

2. *and if*

$$32\gamma^2 \left( \frac{\rho^2 + 1}{(1-\gamma)^2 (1-\rho^2)^2} \left( \frac{\delta_1^2}{\sigma_1^2} + \frac{\delta_2^2}{\sigma_2^2} \right) - 4 \frac{\rho \delta_1 \delta_2}{\sigma_1 \sigma_2 (1-\gamma)^2 (1-\rho^2)^2} + 4 f_2^2(t) \left( \sigma_1^2 + \beta^2 \sigma_2^2 \right) \right.$$
$$\left. + \frac{4 f_2(t)}{(1-\gamma)(1-\rho^2)} \left( -\alpha - \rho \left( \frac{\delta_2 \sigma_1}{\sigma_2} + \frac{\beta \delta_1 \sigma_2}{\sigma_1} \right) \right) \right) (T-t) < \frac{\alpha}{2\sigma_\beta^2},$$

3. *and if*

$$\frac{4\gamma}{(1-\gamma)(1-\rho^2)} \left( \frac{\delta_1^2}{\sigma_1^2} + \frac{\delta_2^2}{\sigma_2^2} - 2\rho \frac{\delta_1 \delta_2}{\sigma_1 \sigma_2} \right) (T-t) - 8\gamma\alpha f_2(t)(T-t) < \frac{\alpha}{2\sigma_\beta^2},$$

*the value function of the optimal stochastic problem is given by*

$$u(t, w, x, y) = \frac{1}{\gamma(1-\gamma)} w^\gamma \exp\{(1-\gamma)(f_2(t) z^2 + f_1(t) z + f_0(t))\},$$

*where the functions $f_0$, $f_1$, and $f_2$ are given in lemma 2.2 and lemma 2.3. Furthermore, the optimal control pair is given by (2.19) and (2.20).*

Since the proof of Theorem 2.1 uses exactly the same sequence of arguments as in Benth and Karlsen (2005), we do not present it. The three conditions on the parameters translate into a limit on the time-horizon. Beyond this bound, the value function may blow up.

Next, we briefly describe the explicit solution for the exponential utility function since the model is slightly more general than in Tourin and Yan (2013).

## 2.1 The solution for the exponential utility function

Now, we allow the wealth variable to be negative and no longer impose constraint (2.8). The terminal condition coupled with the HJB equation (2.11) is $u(T,w,x,y) = -e^{-\gamma w}$, where $\gamma > 0$.

We use the ansatz proposed by Pliska (1986)

$$u(t,w,x,y) = -e^{-\gamma w e^{r(T-t)}} h(t,z),$$

where $h(T,z) = 1$ for all $z$.

**Lemma 2.4.** *1.*

$$h(t,z) = e^{f_2(t)z^2 + f_1(t)z + f_0(t)},$$

*where*

$$\begin{aligned}
f_2(t) &= -\frac{1}{2(1-\rho^2)}\left(\frac{\delta_1^2}{\sigma_1^2} + \frac{\delta_2^2}{\sigma_2^2} - 2\rho\frac{\delta_1\delta_2}{\sigma_1\sigma_2}\right)(T-t), \\
f_1(t) &= -\frac{1}{1-\rho^2}\left[\frac{(\mu_1 - r)\delta_1}{\sigma_1^2} + \frac{(\mu_2 - r)\delta_2}{\sigma_2^2} - \rho\frac{\delta_2(\mu_1 - r) + \delta_1(\mu_2 - r)}{\sigma_1\sigma_2}\right](T-t) \\
&\quad + \frac{1}{1-\rho^2}\left(\frac{\delta_1^2}{\sigma_1^2} + \frac{\delta_2^2}{\sigma_2^2} - 2\rho\frac{\delta_1\delta_2}{\sigma_1\sigma_2}\right)\left(-r(1+\beta) - b + \frac{1}{2}\left(\sigma_1^2 + \beta\sigma_2^2\right)\right)\frac{(T-t)^2}{2}, \\
f_0(t) &= -\frac{1}{2(1-\rho^2)}\left[\frac{(\mu_1 - r)^2}{\sigma_1^2} + \frac{(\mu_2 - r)^2}{\sigma_2^2} - 2\rho\frac{(\mu_1 - r)(\mu_2 - r)}{\sigma_1\sigma_2}\right](T-t) \\
&\quad + \left(r(1+\beta) + b - \frac{1}{2}\left(\sigma_1^2 + \beta\sigma_2^2\right)\right)\int_t^T f_1(s)\,ds \\
&\quad + \left(\sigma_1^2 + \beta^2\sigma_2^2 + 2\beta\rho\sigma_1\sigma_2\right)\int_t^T f_2(s)\,ds.
\end{aligned}$$

*2. Furthermore the controls are of the form*

$$\begin{aligned}
\pi_1^*(t,w,z) &= e^{-r(T-t)}\left[\frac{\mu_1 - r + \delta_1 z}{\sigma_1^2 w\gamma(1-\rho^2)} - \rho\frac{\mu_2 - r + \delta_2 z}{\sigma_1\sigma_2 w\gamma(1-\rho^2)} + \frac{2f_2(t)z + f_1(t)}{w\gamma}\right], \\
\pi_2^*(t,w,z) &= e^{-r(T-t)}\left[\frac{\mu_2 - r + \delta_2 z}{\sigma_2^2 w\gamma(1-\rho^2)} - \rho\frac{\mu_1 - r + \delta_1 z}{\sigma_1\sigma_2 w\gamma(1-\rho^2)} + \beta\frac{2f_2(t)z + f_1(t)}{w\gamma}\right].
\end{aligned}$$

Note that the rate of return is still defined by (2.21). However, the optimal controls are now inversely proportional to the current value of the portfolio. Thus, arbitrarily dividing the initial wealth by 10 causes the excess rate of return to be multiplied by 10. Consequently, the rate of return is a meaningless measure of performance in the case of the exponential utility function since it can be made arbitrarily large.

# 3 A pairs trading model with local volatility

Generally, volatility coefficients can be modelled as deterministic functions of the time, the current stock price and some factors. We refer to Gatheral (2006) for a general introduction to stochastic volatility models. In our framework, we can replace the constant volatility coefficients $\sigma_1, \sigma_2$ by the functions $\sigma_1\left(t, S_t, Y_s^1\right)$ and $\sigma_2\left(t, S_t, Y_s^2\right)$, where for instance, the factors $Y_s^1, Y_s^2$ are mean-reverting. Since, in practice, the computational burden of a general stochastic volatility model is enormous, we focus in this paper on the particular case where the volatility coefficients are respectively modelled by two deterministic functions of the time and current asset log-prices, $\sigma_1\left(t, x_t\right), \sigma_2\left(t, y_t\right)$. A simple choice is provided by the well-known Constant Elasticity of Variance (CEV) model proposed by Cox (1975) and Cox and Ross (1976), namely

$$\sigma_1\left(t, x\right) = \sigma_1 e^{\theta_1 x}, \quad \sigma_2\left(t, y\right) = \sigma_2 e^{\theta_2 y},$$

where $\theta_1 \in (-1, 0)$, $\theta_2 \in (-1, 0)$, $\sigma_1 > 0$, and $\sigma_2 > 0$. There is a large amount of literature on CEV models and we do not provide the reader with an exhaustive list of articles. We simply refer, among many others, to Beckers (1980), Delbaen and Shirakawa (2002), and Geman and Shih (2009) for an overview of the mathematical, modelling and econometric aspects of this subject. Unfortunately, it is not a realistic model for the dynamics of the volatility but, at least, it captures the negative link between volatility and stock price observed on real financial data. The combination of the Merton problem with a CEV volatility model was solved explicitly by Gao (2009). However, for the pairs trading strategies under the CEV model, we are unfortunately unable to derive a fully explicit solution.

Next, we present the model and derive the linear PDE characterizing the solution. The dynamics of the state variables $(W_s, x_s, y_s)$ are now given by

$$\begin{aligned} dW_s &= \pi_s^1\left(\mu_1 + \delta_1 z_s\right) W_s ds + \pi_s^2\left(\mu_2 + \delta_2 z_s\right) W_s ds \\ &\quad + \pi_s^1 \sigma_1\left(s, x_s\right) W_s dB_s^1 + \pi_s^2 \sigma_2\left(s, y_s\right) W_s dB_s^2, \\ dx_s &= \left(\mu_1 - \frac{\sigma_1^2\left(s, x_s\right)}{2} + \delta_1 z_s\right) ds + \sigma_1\left(s, x_s\right) dB_s^1, \\ dy_s &= \left(\mu_2 - \frac{\sigma_2^2\left(s, y_s\right)}{2} + \delta_2 z_s\right) ds + \sigma_2\left(s, y_s\right) dB_s^2, \\ W_t &= w,\ x_t = x,\ y_t = y, \end{aligned}$$

where $z$ is defined in (2.3). We note that $z$ satisfies now the dynamics

$$\begin{aligned} dz_t = &\left(b + \mu_1 + \beta\mu_2 - \frac{1}{2}\sigma_1^2\left(t, x_t\right) - \frac{1}{2}\beta\sigma_2^2\left(t, y_t\right) + \delta_1 z_t + \beta\delta_2 z_t\right) dt \\ &+ \sigma_1\left(t, x_t\right) dB_t^1 + \beta\sigma_2\left(t, y_t\right) dB_t^2. \end{aligned}$$

Furthermore, we notice that the auxiliary process $z_s'$, defined as

$$z_s' = z_s + \int_t^s e^{\alpha(u-s)}\left(\frac{1}{2}\sigma_1^2\left(u, x_u\right) + \frac{1}{2}\beta\sigma_2^2\left(u, y_u\right)\right) du, \tag{3.1}$$

is mean-reverting and that its long-run mean is the constant $\frac{1}{\alpha}\left(b + \mu_1 + \beta\mu_2\right)$.

As in the constant volatility case, we require the wealth process to satisfy (2.8). The value function $u$ is defined in (2.9), where the utility function $U$ is given in (2.10).

We then rewrite $u$ as $u(t,w,x,y) = \frac{1}{\gamma} w^{\gamma} g(t,x,y)$, where $g$ satisfies the PDE

$$\begin{aligned} &- g_t - \sup_{\pi_1,\pi_2} \Big\{ [\pi_1(\mu_1 + \delta_1 z) + \pi_2(\mu_2 + \delta_2 z)]\gamma g + r(1 - \pi_1 - \pi_2)\gamma g \\ &+ \left(\mu_1 - \frac{1}{2}\sigma_1^2(t,x) + \delta_1 z\right) g_x + \left(\mu_2 - \frac{1}{2}\sigma_2^2(t,x) + \delta_2 z\right) g_y \\ &+ \left(\pi_1\sigma_1^2(t,x) + \pi_2\rho\sigma_1(t,x)\sigma_2(t,x)\right)\gamma g_x + (\pi_2\sigma_2^2(t,x) + \pi_1\rho\sigma_1(t,x)\sigma_2(t,x))\gamma g_y \\ &+ \frac{1}{2}\sigma_1^2(t,x) g_{xx} + \frac{1}{2}\sigma_2^2(t,x) g_{yy} + \rho\sigma_1(t,x)\sigma_2(t,x) g_{xy} \\ &+ \frac{1}{2}\left(\pi_1^2\sigma_1^2(t,x) + \pi_2^2\sigma_2^2(t,x) + 2\pi_1\pi_2\rho\sigma_1(t,x)\sigma_2(t,x)\right)\gamma(\gamma-1) g \Big\} = 0, \end{aligned} \tag{3.2}$$

for all $0 \le t < T$, $x$, and $y$, coupled with $g(T,x,y) = 1$, for all $x$ and $y$.

We compute the optimal controls by solving the optimization problem in (3.2), in terms of $g$ and its partial derivatives

$$\pi_1^* = \frac{\mu_1 - r + \delta_1 z}{\sigma_1^2(t,x)(1-\gamma)(1-\rho^2)} - \rho\frac{\mu_2 - r + \delta_2 z}{\sigma_1(t,x)\sigma_2(t,x)(1-\gamma)(1-\rho^2)} + \frac{g_x}{(1-\gamma)g}, \tag{3.3}$$

$$\pi_2^* = \frac{(\mu_2 - r + \delta_2 z)}{\sigma_2^2(t,x)(1-\gamma)(1-\rho^2)} - \rho\frac{\mu_1 - r + \delta_1 z}{\sigma_1(t,x)\sigma_2(t,x)(1-\gamma)(1-\rho^2)} + \frac{g_y}{(1-\gamma)g}. \tag{3.4}$$

Substituting the controls (3.3) and (3.4) into (3.2) yields

$$\begin{aligned} &- g_t + \frac{\gamma}{2(\gamma-1)(1-\rho^2)} \\ &\left[\frac{(\mu_1 + \delta_1 z - r)^2}{\sigma_1^2(t,x)} + \frac{(\mu_2 + \delta_2 z - r)^2}{\sigma_2^2(t,x)} - 2\rho\frac{(\mu_1 - r + \delta_1 z)(\mu_2 - r + \delta_2 z)}{\sigma_1\sigma_2}\right] g \\ &- r\gamma g \\ &+ \frac{(\mu_1 + \delta_1 z)}{\gamma - 1} g_x - \frac{r\gamma}{\gamma-1} g_x + \frac{1}{2}\sigma_1^2(t,x) g_x \\ &+ \frac{(\mu_2 + \delta_2 z)}{\gamma - 1} g_y - \frac{r\gamma}{\gamma-1} g_y + \frac{1}{2}\sigma_2^2(t,x) g_y \\ &+ \frac{\gamma}{2(\gamma-1)}\left(\sigma_1^2(t,x)\frac{g_x^2}{g} + \sigma_2^2(t,x)\frac{g_y^2}{g} + \rho\sigma_1(t,x)\sigma_2(t,x)\frac{g_x g_y}{g}\right) \\ &- \frac{1}{2}\sigma_1^2(t,x) g_{xx} - \frac{1}{2}\sigma_2^2(t,y) g_{yy} - \rho\sigma_1(t,x)\sigma_2(t,x) g_{xy} = 0. \end{aligned} \tag{3.5}$$

We remove the non-linearities in (3.5) by using the same technique as in the previous section; we define for $\tau = T - t$, and $x, y$, the function $\phi(\tau,x,y)$, such that $g(t,x,y) = \frac{1}{1-\gamma}\phi^{1-\gamma}(\tau,x,y)$. The linear PDE satisfied by $\phi$ reads

$$\begin{aligned} &\phi_\tau - \frac{\gamma}{2(\gamma-1)^2(1-\rho^2)} \\ &\left[\frac{(\mu_1 + \delta_1 z - r)^2}{\sigma_1^2(t,x)} + \frac{(\mu_2 + \delta_2 z - r)^2}{\sigma_2^2(t,x)} - 2\rho\frac{(\mu_1 - r + \delta_1 z)(\mu_2 - r + \delta_2 z)}{\sigma_1(t,x)\sigma_2(t,x)}\right]\phi \\ &- \frac{r\gamma}{1-\gamma}\phi + \frac{\mu_1 + \delta_1 z}{\gamma-1}\phi_x - \frac{r\gamma}{\gamma-1}\phi_x + \frac{1}{2}\sigma_1^2(t,x)\phi_x \\ &+ \frac{\mu_2 + \delta_2 z}{\gamma-1}\phi_y - \frac{r\gamma}{\gamma-1}\phi_y + \frac{1}{2}\sigma_2^2(t,y)\phi_y \end{aligned}$$

$$-\frac{1}{2}\sigma_1^2(t,x)\phi_{xx}-\frac{1}{2}\sigma_2^2(t,x)\phi_{yy}-\rho\sigma_1(t,x)\sigma_2(t,x)\phi_{xy}=0.$$

In particular, for the choice $\sigma_1(t,x)=\sigma_1 e^{\theta_1 x}$ and $\sigma_2(t,x)=\sigma_2 e^{\theta_2 y}$, it becomes

$$\begin{aligned}&\phi_\tau-\frac{\gamma}{2(\gamma-1)^2(1-\rho^2)}\\&\left[\frac{(\mu_1+\delta_1 z-r)^2}{\sigma_1^2 e^{2\theta_1 x}}+\frac{(\mu_2+\delta_2 z-r)^2}{\sigma_2^2 e^{2\theta_2 y}}-2\rho\frac{(\mu_1-r+\delta_1 z)(\mu_2-r+\delta_2 z)}{\sigma_1(t,x)\sigma_2(t,x)}\right]\phi\\&-\frac{r\gamma}{1-\gamma}\phi+\frac{\mu_1+\delta_1 z}{\gamma-1}\phi_x-\frac{r\gamma}{\gamma-1}\phi_x+\frac{1}{2}\sigma_1^2 e^{2\theta_1 x}\phi_x\\&+\frac{\mu_2+\delta_2 z}{\gamma-1}\phi_y-\frac{r\gamma}{\gamma-1}\phi_y+\frac{1}{2}\sigma_2^2 e^{2\theta_2 y}\phi_y-\\&\frac{1}{2}\sigma_1^2 e^{2\theta_1 x}\phi_{xx}-\frac{1}{2}\sigma_2^2 e^{2\theta_2 y}\phi_{yy}-\rho\sigma_1\sigma_2 e^{\theta_1 x+\theta_2 y}\phi_{xy}=0,\end{aligned}\tag{3.6}$$

and must be coupled with the initial condition $\phi(0,x,y)=(1-\gamma)^{\frac{1}{1-\gamma}}$.

Finally, we rewrite the controls in terms of the function $\phi$ and its partial derivatives:

$$\pi_1^*=\frac{\mu_1-r+\delta_1 z}{\sigma_1^2(t,x)(1-\gamma)(1-\rho^2)}-\rho\frac{\mu_2-r+\delta_2 z}{\sigma_1(t,x)\sigma_2(t,x)(1-\gamma)(1-\rho^2)}+\frac{\phi_x}{\phi},\tag{3.7}$$

$$\pi_2^*=\frac{\mu_2-r+\delta_2 z}{\sigma_2^2(t,x)(1-\gamma)(1-\rho^2)}-\rho\frac{\mu_1-r+\delta_1 z}{\sigma_1(t,x)\sigma_2(t,x)(1-\gamma)(1-\rho^2)}+\frac{\phi_y}{\phi}.\tag{3.8}$$

We see that the first two components in each control are still known in closed form and are essentially the same as in the constant volatility case, except that the volatility coefficients are no longer constant but depend on the log prices and the time-variable. However, the third term is no longer known in closed form and computing it requires approximating the above PDE numerically. As in the constant volatility case, we expect the computed solution to blow up as the time horizon becomes large, and we observe this phenomenon experimentally. Unfortunately, when the local volatility is no longer constant, we are unable to derive in closed form a bound on the time horizon. In the next section, we present our numerical method and the tests we conduct.

# 4 Numerical experiments

Since the PDE (3.6) cannot be solved analytically, we use a purely implicit Finite Difference (FD in short) scheme to approximate the value function and optimal trading strategy (3.7) and (3.8). Most of the finite differences are chosen not only consistent but also unconditionally monotone as in Barles and Souganidis (1991), except for the cross-partial derivative, which is simply approximated on a 7-point stencil as in Ma and Forsyth (2014). We then verify the stability of the scheme through successive refinements of the mesh. We also impose zero Neumann boundary conditions at the edges of the computational domain. These values underestimate greatly the actual growth of the solution for large values of $z$. Indeed, we expect $\phi$ to grow exponentially in the variable $z^2$ as $z\to\pm\infty$. Consequently, the numerical error near the boundary is substantially greater than inside the domain and we must take this into account when setting the size of the computational domain. This type of numerical scheme is quite standard and we refer to the appendix for its detailed description.

The algorithm performs well, which is most likely due to the linearity of the degenerate diffusion operator. We start by performing a collection of experiments for an arbitrary set of parameters before turning to market data. In practice, we set

$$\beta = -0.6,\ a = -0.01,\ b = -0.01,\ \sigma_1 = 0.3,\ \sigma_2 = 0.35,\ T = 1 \text{ year},\ r = 0.01,$$
$$\delta_1 = -0.1,\ \delta_2 = 0.1,\ \mu_1 = 0.2,\ \mu_2 = 0.08,\ \gamma = 0.1,\ \rho = 0.5.$$

We also work in the finite domain $[1,5] \times [1,5]$ for the spatial variables.

First of all, to ensure that our algorithm works, we compare, for the above set of parameters, our solution with the closed form expression that is available in the constant volatility case, i.e. for $\theta_1 = \theta_2 = 0$, and compute the error. For the discretization steps $\Delta t = 1/251$, $\Delta x = \Delta y = 0.1$, we find a $L^\infty$ error equal to 0.0595 and we verify that, as expected, the maximum error occurs on the edge of the domain.

Secondly, we compute the value function and the optimal control, for the CEV model with $\theta_1 = -0.2$, $\theta_2 = -0.15$, $\sigma_1 = 0.4$, $\sigma_2 = 0.45$. We leave the other parameters, the computational domain, and the discretization steps unchanged. We show the optimal controls at $t = 0$ as functions of $x$ and $y$ in Figure 4.1-4.2, and the function $\phi$ in Figure 4.3.

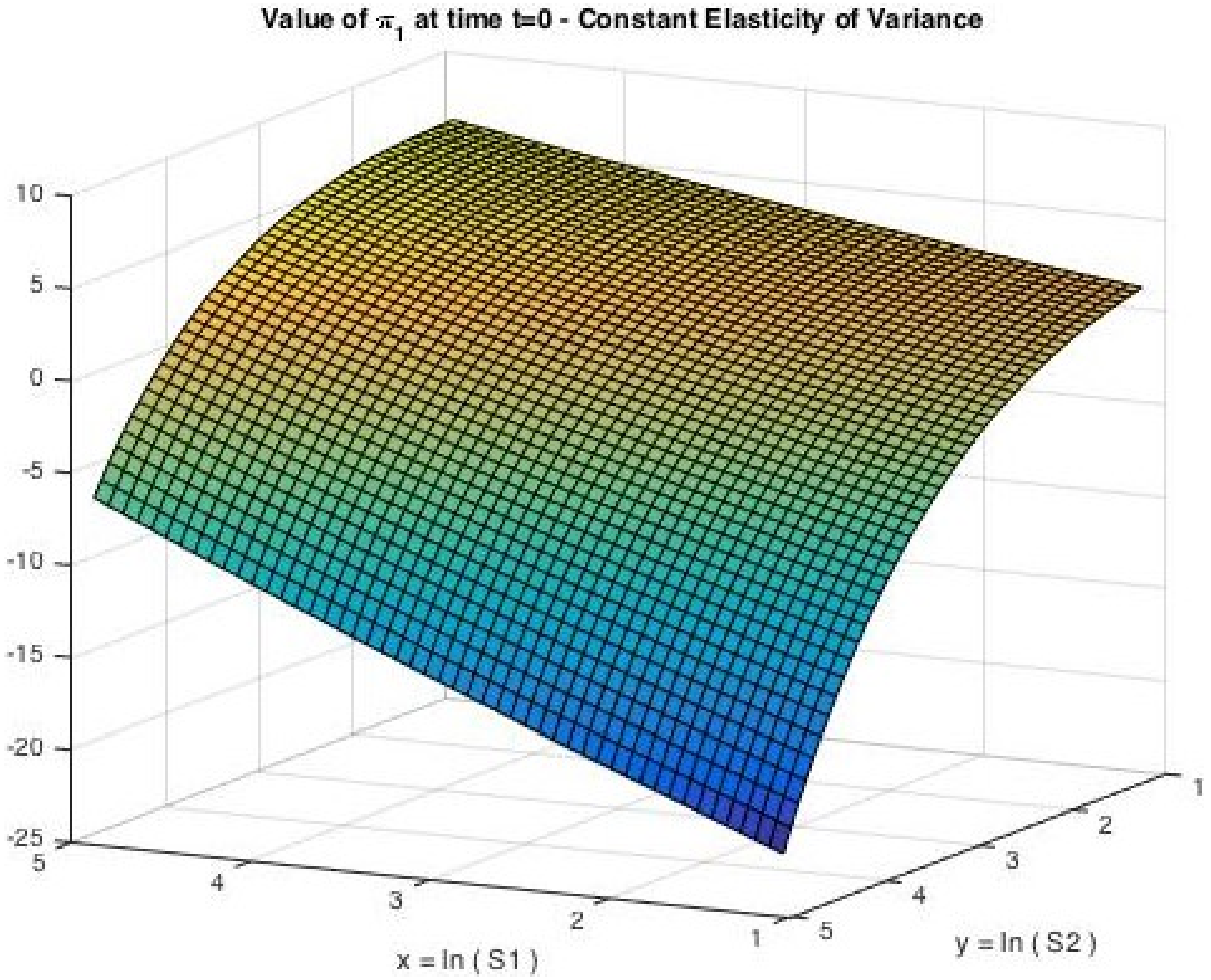

Figure 4.1: This shows the fraction of wealth invested in the first risky asset for a time horizon of one year, as a function of the assets' log-prices, for the CEV model. It was computed by the FD scheme on a $50 \times 50$ spatial mesh and with 251 time steps. The parameters are $r = 0.01$, $\mu_1 = 0.2$, $\mu_2 = 0.08$, $\sigma_1 = 0.4$, $\sigma_2 = 0.45$, $\beta = -0.6$, $a = -0.01$, $b = -0.01$, $\delta_1 = -0.1$, $\delta_2 = 0.1$, $\theta_1 = -0.2$, $\theta_2 = -0.15$, $\gamma = 0.1$, $\rho = 0.5$

Thirdly, we illustrate the applicability of the model with historical data. Our first experiment is conducted with one year of daily data, from January 2nd, 2013 to December 31st, 2013 on two assets, the Market Vectors Gold Miners ETF, with ticker symbol GDX and the SPDR Gold Shares

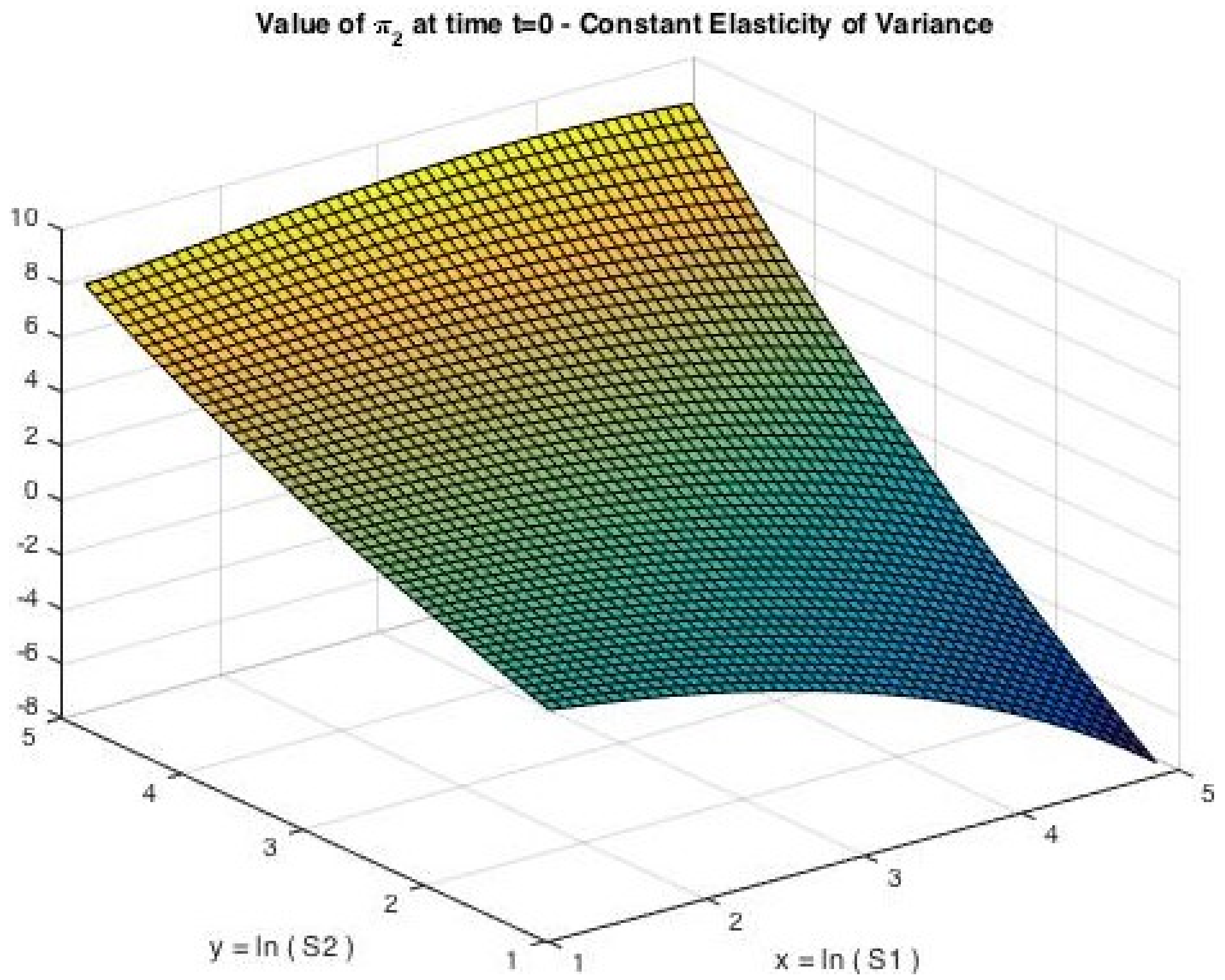

Figure 4.2: This shows the fraction of wealth invested in the second risky asset for a time horizon of one year, as a function of the assets' log-prices for the CEV model. It was computed by the FD scheme on a $50 \times 50$ spatial mesh and with 251 time steps. The parameters are $r = 0.01$, $\mu_1 = 0.2$, $\mu_2 = 0.08$, $\sigma_1 = 0.4$, $\sigma_2 = 0.45$, $\beta = -0.6$, $a = -0.01$, $b = -0.01$, $\delta_1 = -0.1$, $\delta_2 = 0.1$, $\theta_1 = -0.2$, $\theta_2 = -0.15$, $\gamma = 0.1$, $\rho = 0.5$.

ETF with ticker symbol GLD. To test for co-integration, we use on the one hand the two-step Engle-Granger approach combined with both the Augmented Dickey-Fuller and Phillips-Perron tests, and on the other hand, the Johansen co-integration test. The data pass all these tests at the 5% significance level. We want to mention that the use of these tests is somewhat inconsistent with our assumption that the volatility is actually of the CEV type, which is non-stationary. In future work, we are planning to implement the robust Bayesian unit root test proposed by Zhang et al. (2013). Next, we estimate the parameters in the co-integrated asset price model (2.1), (2.2), and (2.3) by the Generalized Method of Moments (GMM in short). We refer to the Appendix for a detailed description of our methodology. We obtain the following annualized parameters

$$\mu_1 = 0.55,\ \sigma_1 = 0.76,\ \delta_1 = -2.57,\ \theta_1 = -0.16,\ \mu_2 = 0.29,\ \sigma_2 = 0.63,$$
$$\delta_2 = -1.25,\ \theta_2 = -0.22,\ a = -0.39,\ b = 0.014,\ \beta = -0.51,\ \rho = 0.78.$$

The objective value obtained by the GMM in the just-identified case at the optimum is $J_T(\lambda^*) = 3.73.10^{-8}$ and the process $z'$ defined in (3.1) passes Augmented Dickey-Fuller stationarity test. In the over-identified case, we find $J_T(\lambda^*) = 0.035$. The computed strategies generate a profit of roughly 2000%. We show the prices in Figure 4.4, the computed trading strategies in Figure 4.5, and the cumulative Profit and Loss in Figure 4.6. Note that both optimal positions are short for this this data, due to the fact that it has a significant downward trend over this time period.

Finally, we perform a simple out-of-sample test, inspired from the preliminary experiments conducted for the constant volatility case and the exponential utility function in Lee (2015) on cross-listed Chinese stocks. Specifically, we collect daily data for the Bank of China shares listed respectively on the New York Stock Exchange (NYSE in short) and the Hong Kong Stock Exchange (HKSE in short), from July 08, 2010 to July 07, 2015. The ticker symbols are respectively OTC:BACHY and HKSE:3988. Since the Hong Kong Dollar has been pegged to the US Dollar since September 1983, resulting in an approximately constant exhange rate between these two currencies, we do not convert the share prices quoted on the Hong Kong stock exchange into US Dollars. For the purpose of conducting an out-of-sample test, we split the whole time period into two subintervals; the first one, spanning four years, from July 08, 2010 to July 07, 2014, is used for estimating the parameters in the model by the GMM. Then, we apply to the second time interval the computed trading strategies associated with the set of parameters estimated from the data in the first interval.

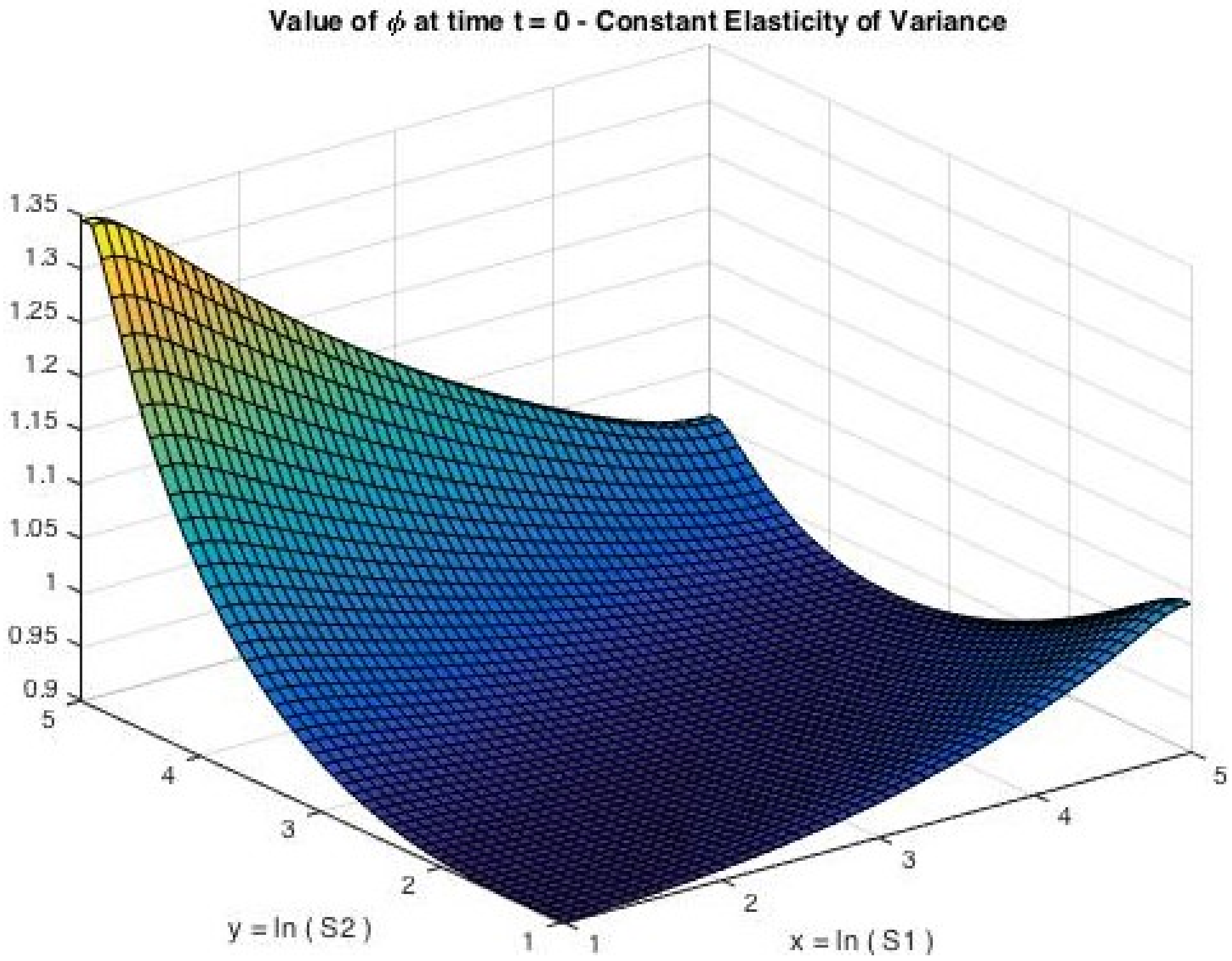

Figure 4.3: This shows the function $\phi$ for a time horizon of one year, as a function of the assets' log-prices for the CEV model. It was computed by the FD scheme on a $50 \times 50$ spatial mesh and with 251 time steps. The values of the parameters are $r = 0.01$, $\mu_1 = 0.2$, $\mu_2 = 0.08$, $\sigma_1 = 0.4$, $\sigma_2 = 0.45$, $\beta = -0.6$, $a = -0.01$, $b = -0.01$, $\delta_1 = -0.1$, $\delta_2 = 0.1$, $\theta_1 = -0.2$, $\theta_2 = -0.15$, $\gamma = 0.1$, $\rho = 0.5$.

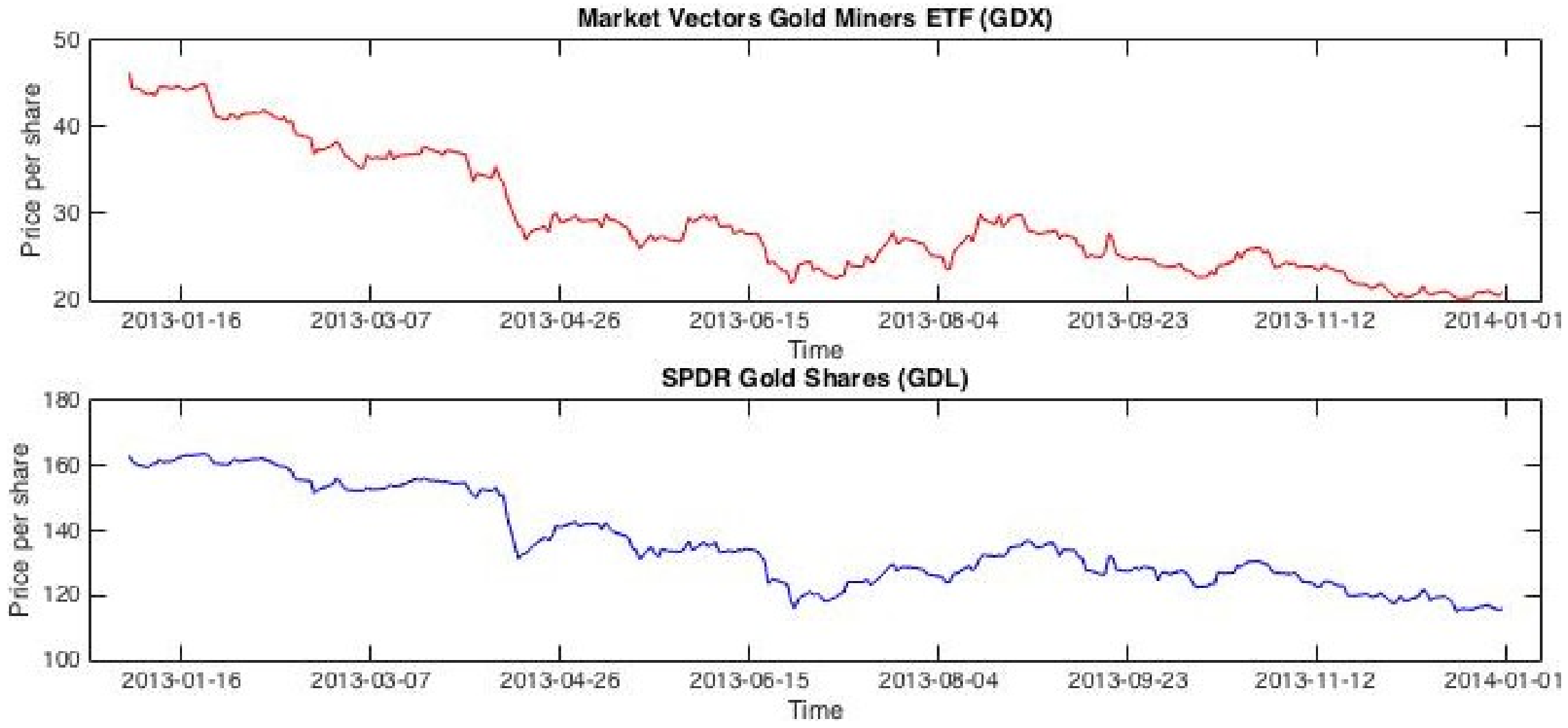

Figure 4.4: Prices of GDX and GLD in 2013 (daily frequency)

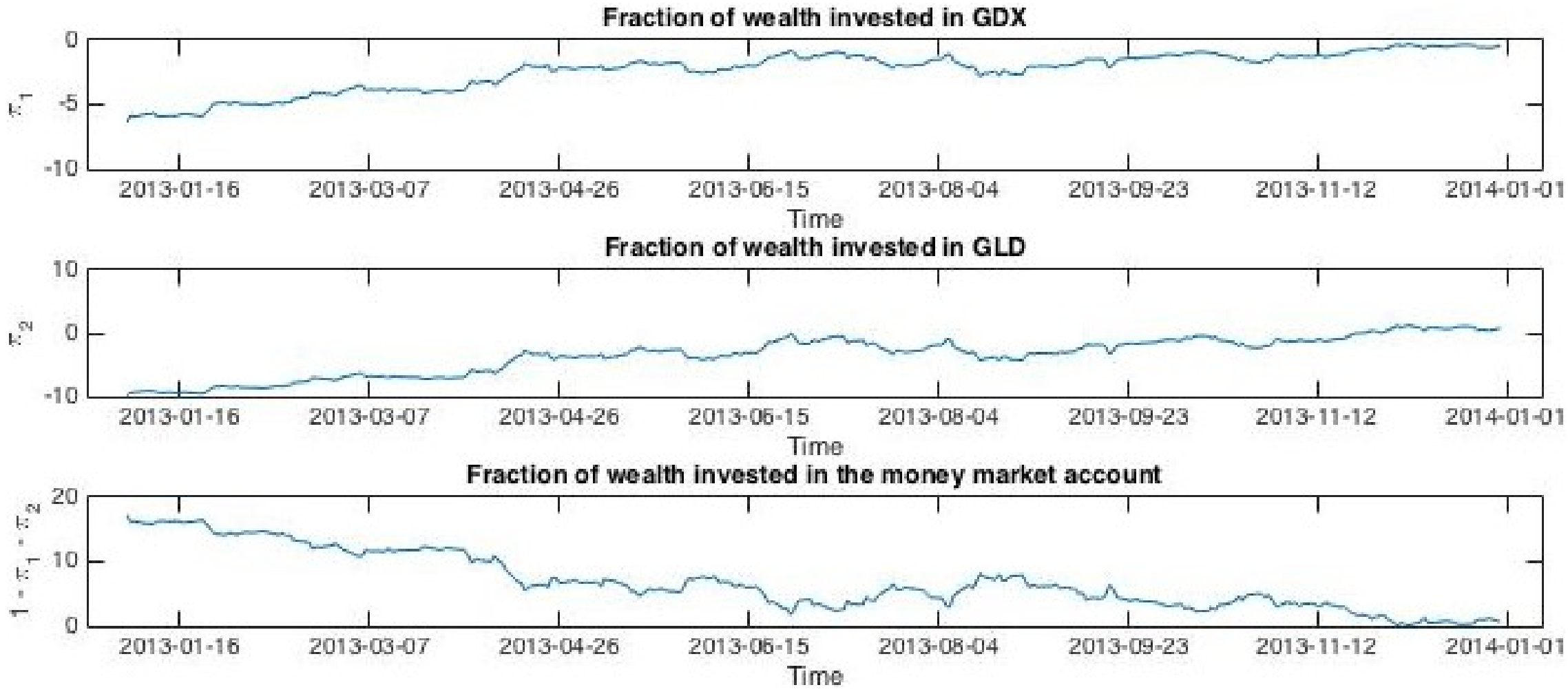

Figure 4.5: This shows the trading strategies computed by the FD scheme with $\Delta t = 1/251$, $\Delta x = 1/60$, $\Delta y = 1/75$, for the asset prices GDX and GLD. The estimated parameters are $\mu_1 = 0.55$, $\sigma_1 = 0.76$, $\delta_1 = -2.57$, $\theta_1 = -0.16$, $\mu_2 = 0.29$, $\sigma_2 = 0.63$, $\delta_2 = -1.25$, $\theta_2 = -0.22$, $a = -0.39$, $b = 0.014$, $\beta = -0.51$, $\rho = 0.78$.

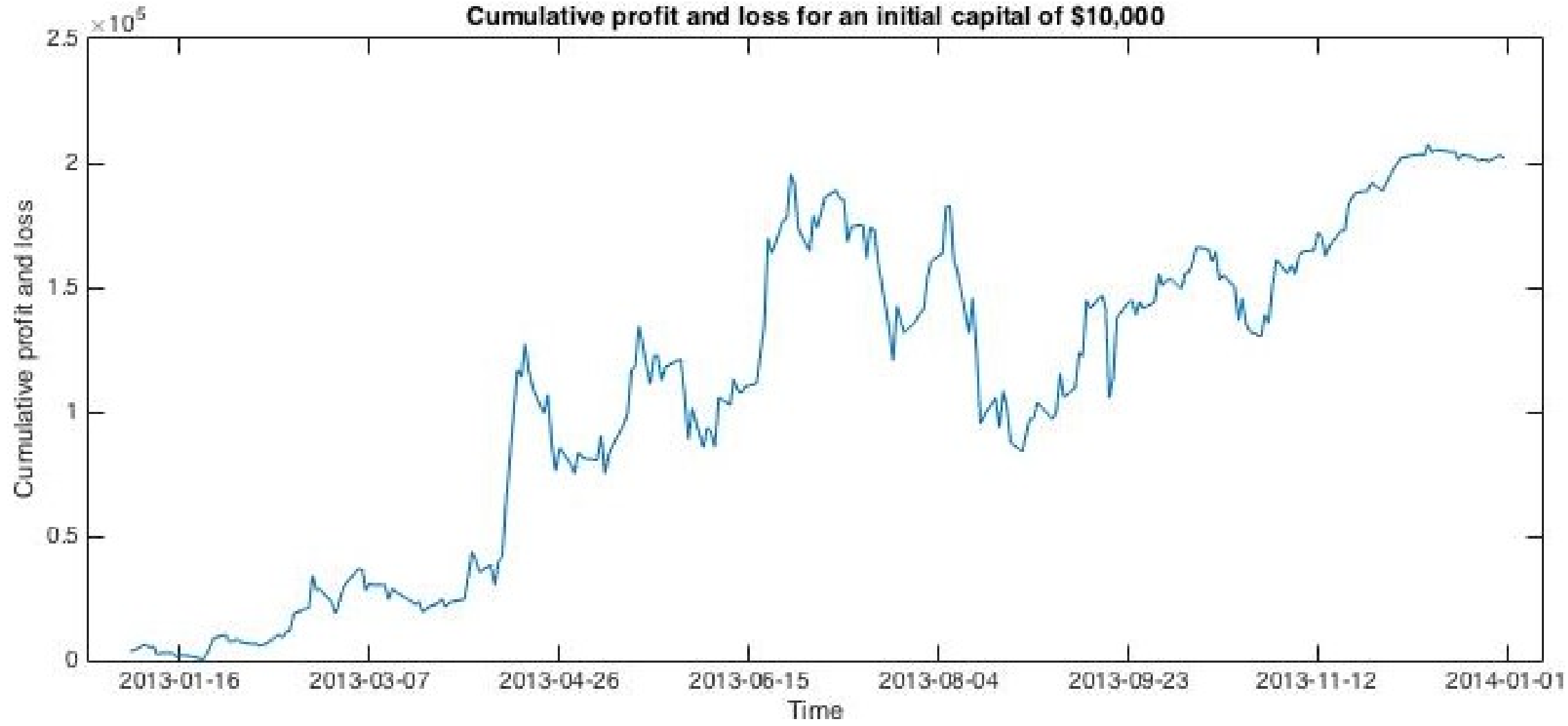

Figure 4.6: This shows the cumulative profit and loss for the above trading strategies for the asset prices GDX and GLD. The estimated parameters are $\mu_1 = 0.55$, $\sigma_1 = 0.76$, $\delta_1 = -2.57$, $\theta_1 = -0.16$, $\mu_2 = 0.29$, $\sigma_2 = 0.63$, $\delta_2 = -1.25$, $\theta_2 = -0.22$, $a = -0.39$, $b = 0.014$, $\beta = -0.51$, $\rho = 0.78$.

The data in each subinterval pass the Augmented Dickey-Fuller, Phillips-Perron and Johansen co-integration tests. The GMM in the just-identified case yields the following results: $\mu_1 = 0.27$,

$\sigma_1 = 0.76$, $\delta_1 = -0.36$, $\theta_1 = -0.16$, $\mu_2 = 0.26$, $\sigma_2 = 0.7$, $\delta_2 = 0.85$, $\theta_2 = -0.37$, $a = -1$, $b = 0.068$, $\beta = -1.12$, $\rho = 0.64$, $J_T(\lambda^*) = 3.06.10^{-08}$. Moreover, the process $z'$ defined in 3.1 passes the Augmented Dickey-Fuller stationarity test. The objective value of the GMM in the over-identified case is $J_T(\lambda^*) = 0.113$. We show the prices in Figure4.7-4.8. In Figure 4.9-4.10, we apply the computed trading strategies to the in-sample data, yielding a profit of roughly 140%. Finally, we present the test on the out-of-sample data in Figure 4.11-4.12. As expected, the profit is lower, roughly equal to 80%. Note that both positions are long, taking advantage of the strong upward trend over this time period.

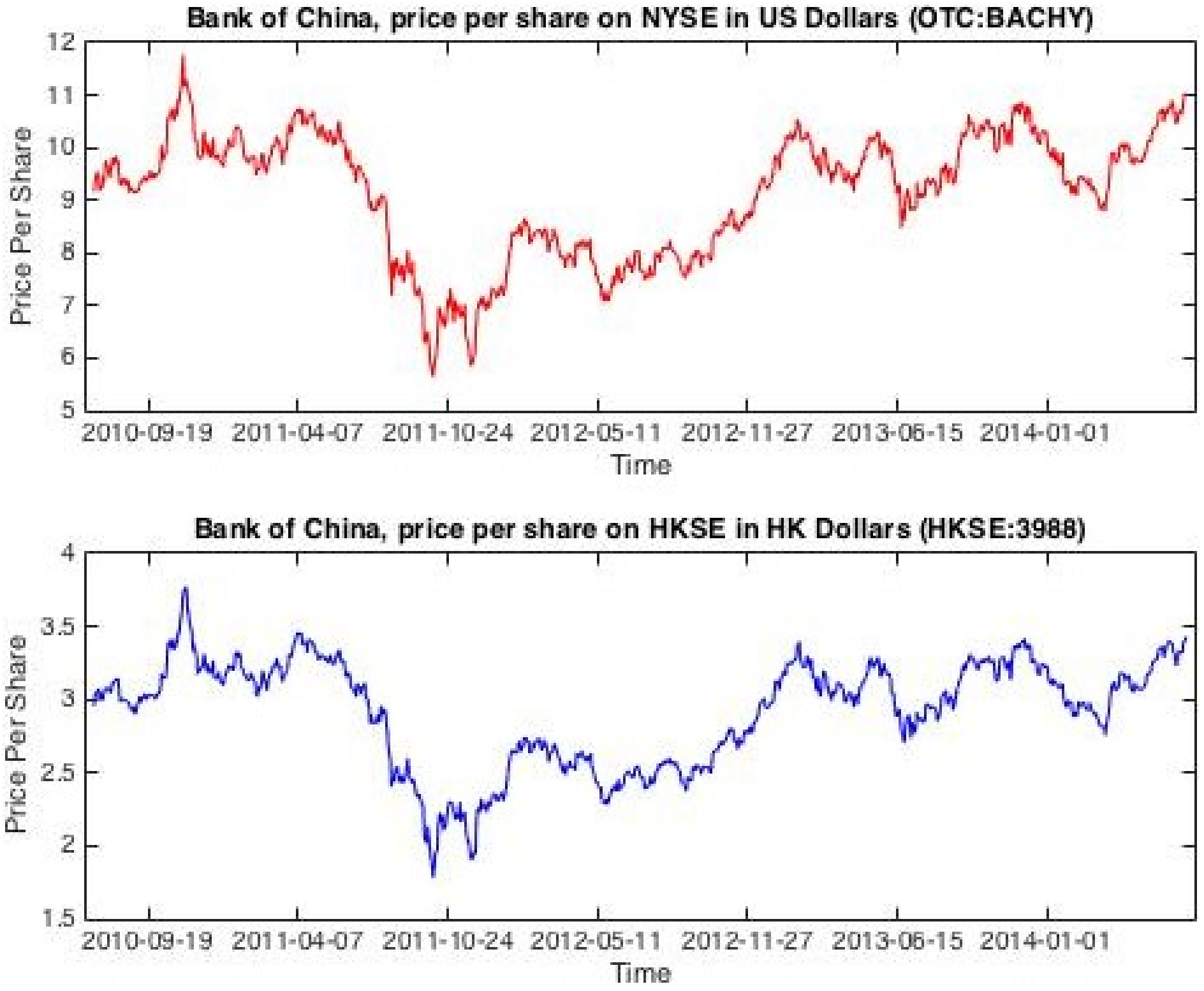

Figure 4.7: This shows the prices for Bank of China on the NYSE and HKSE from July 08, 2010 to July 7, 2014 (in-sample time interval).

## 5 Conclusion

This work constitutes a first step toward incorporating a time-varying volatility model into the pairs trading stochastic control problem proposed in Tourin and Yan (2013). We use here a CEV volatility model and compute numerically the optimal trading strategies. We also illustrate our approach with in-sample and out-of-sample tests. There is clearly a need for an extensive performance study involving a larger number of pairs. Estimated transaction costs should also be subtracted from the cumulative profit and loss.

## Appendix A The Finite Difference scheme

We describe here the fully implicit monotone scheme that we implement. We construct a monotone scheme by using the appropriate forward or backward finite differences for the first

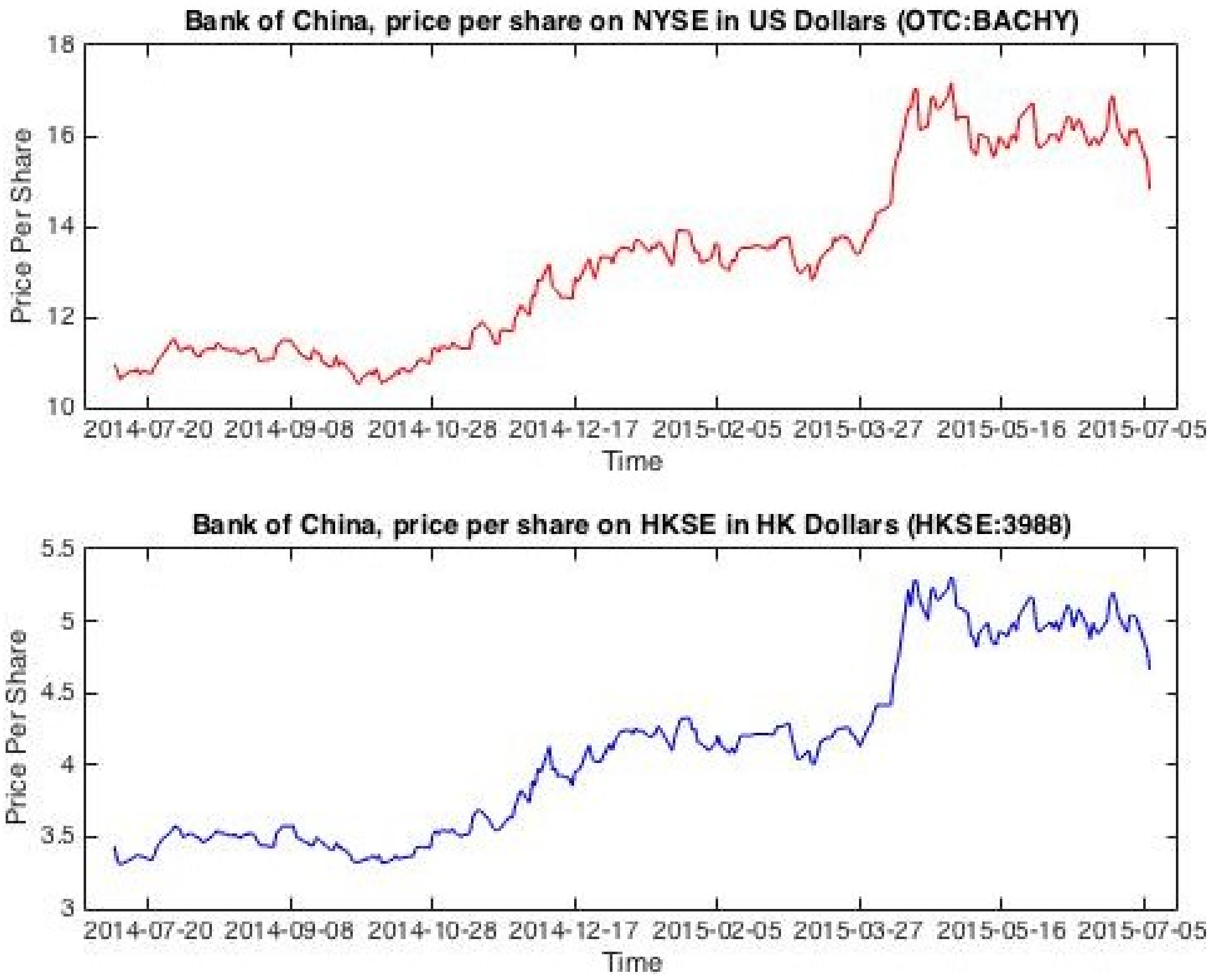

Figure 4.8: This shows the prices for Bank of China on the NYSE and HKSE from July 08, 2014 to July 7, 2015 (out-of-sample time interval).

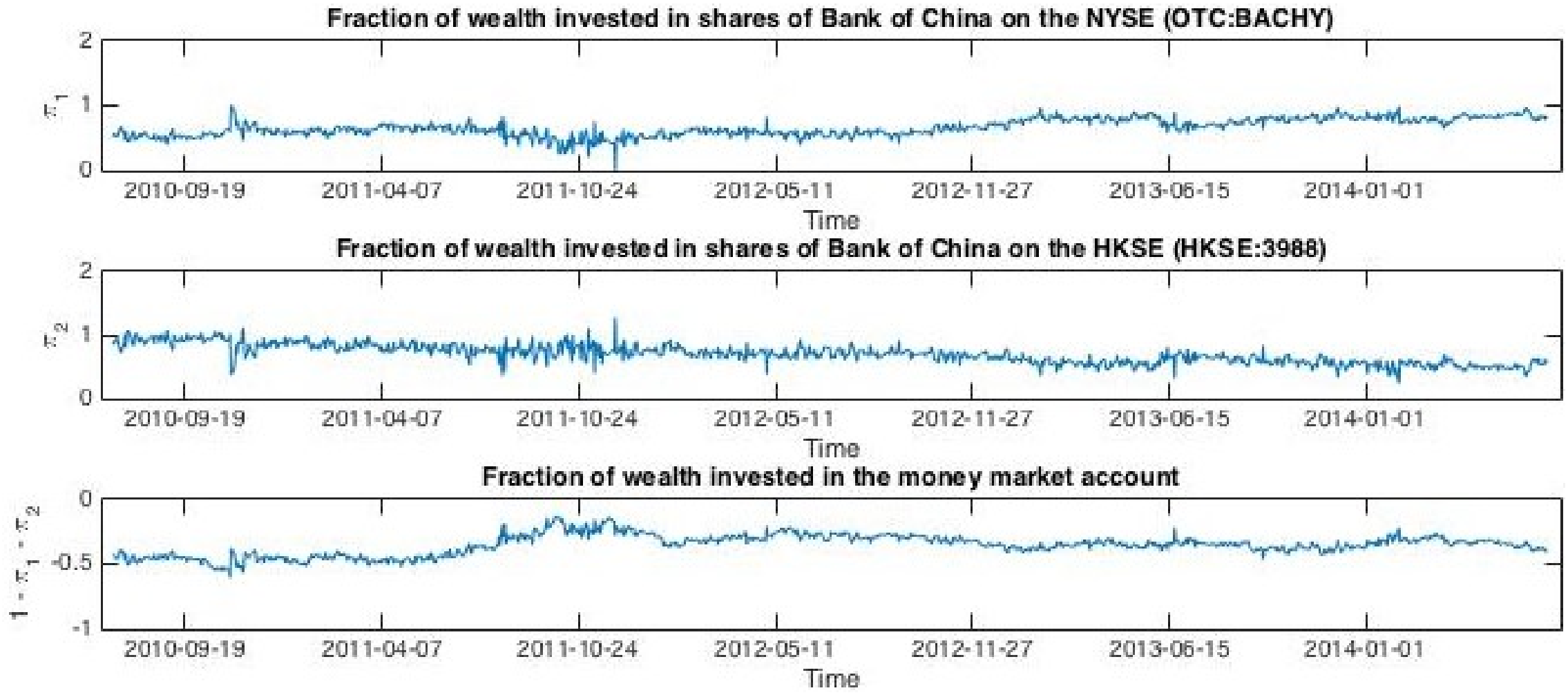

Figure 4.9: This shows the trading strategies computed by the FD scheme for the in-sample test on Bank of China shares, from July 08, 2010 to July 07, 2014. The estimated parameters are $\mu_1 = 0.27$, $\sigma_1 = 0.76$, $\delta_1 = -0.36$, $\theta_1 = -0.16$, $\mu_2 = 0.26$, $\sigma_2 = 0.7$, $\delta_2 = 0.85$, $\theta_2 = -0.37$, $a = -1$, $b = 0.068$, $\beta = -1.12$, $\rho = 0.64$.

spatial derivatives.

First of all, we compute the solution in the bounded domain $[x_{\min},\ x_{\max}] \times [y_{\min},\ y_{\max}]$. Next, let $\Delta t$, $\Delta x$, $\Delta y$ be the discretization steps. We then set

$$x_i = x_{\min} + i\Delta x,$$

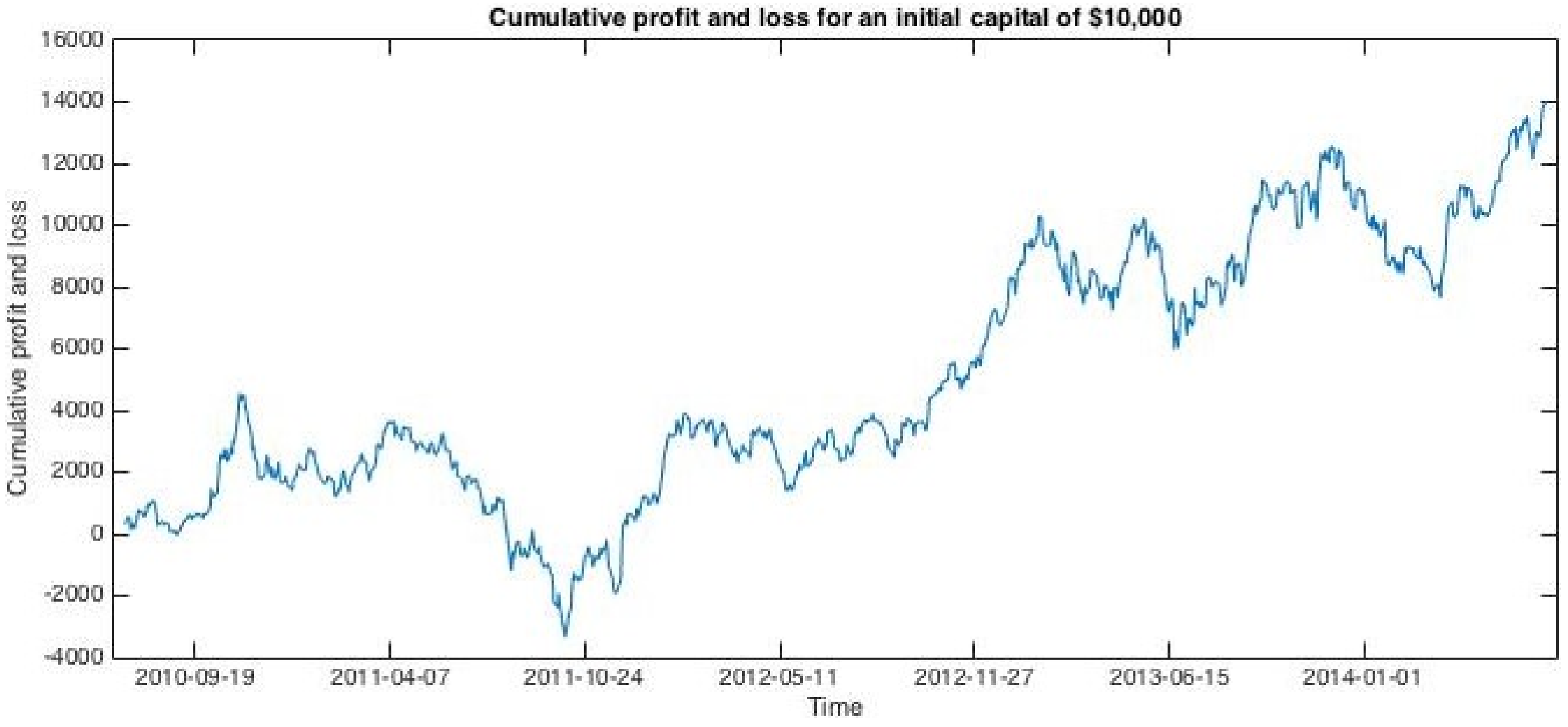

Figure 4.10: This shows the cumulative profit and loss for the in-sample test on Bank of China shares, from July 08, 2010 to July 07, 2014. The estimated parameters are $\mu_1 = 0.27$, $\sigma_1 = 0.76$, $\delta_1 = -0.36$, $\theta_1 = -0.16$, $\mu_2 = 0.26$, $\sigma_2 = 0.7$, $\delta_2 = 0.85$, $\theta_2 = -0.37$, $a = -1$, $b = 0.068$, $\beta = -1.12$, $\rho = 0.64$..

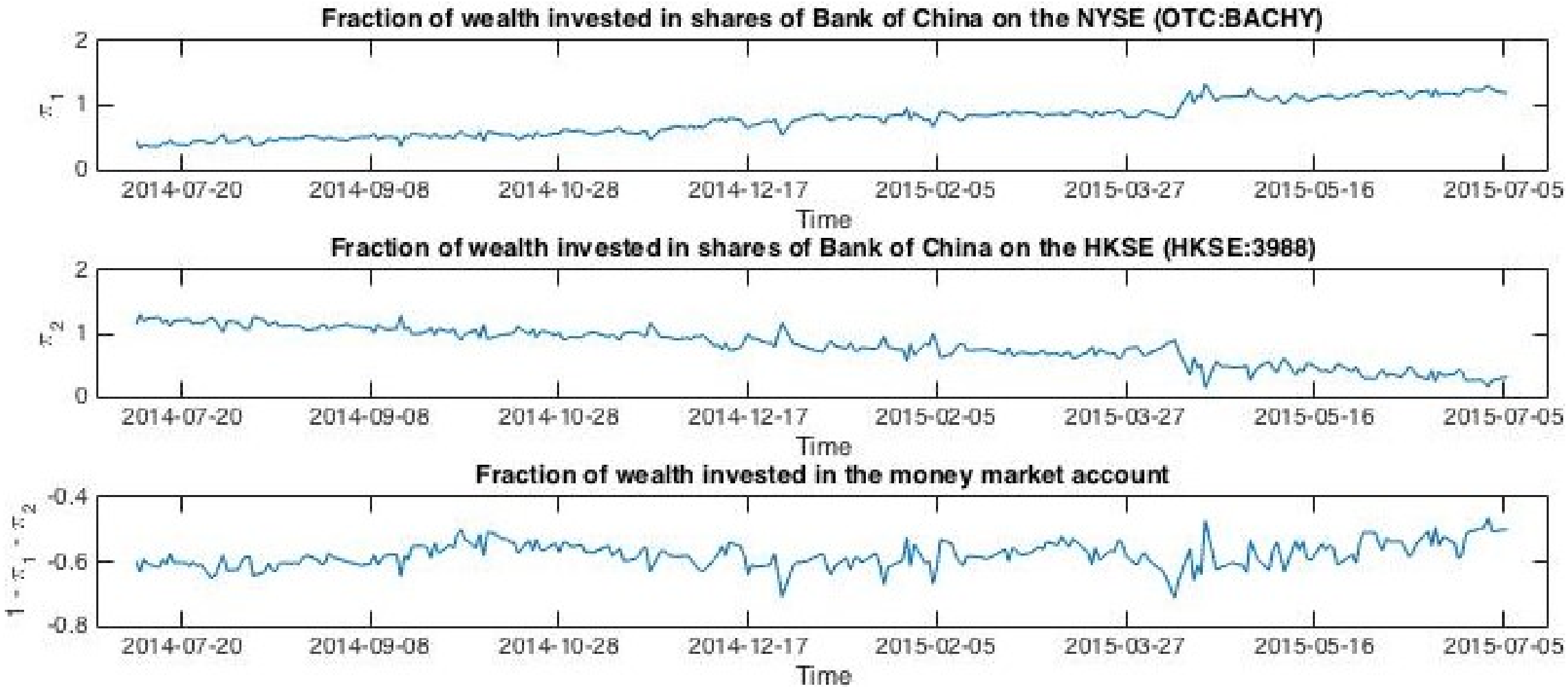

Figure 4.11: This shows the trading strategies computed by the FD scheme for the out-of-sample test on Bank of China shares, from July 08, 2014 to July 07, 2015. The parameters are estimated from the in-sample interval: $\mu_1 = 0.27$, $\sigma_1 = 0.76$, $\delta_1 = -0.36$, $\theta_1 = -0.16$, $\mu_2 = 0.26$, $\sigma_2 = 0.7$, $\delta_2 = 0.85$, $\theta_2 = -0.37$, a$= -1$, $b = 0.068$, $\beta = -1.12$, $\rho = 0.64$.

$$
\begin{aligned}
y_j &= y_{\min} + j\Delta y, \\
z_{i,j}^k &= a + b\,(T - k\Delta t) + x_i + \beta y_j,
\end{aligned}
$$

for $i,\ j \in \{0, \cdots, I\} \times \{0, \cdots, J\}$, with $I\Delta x = x_{\max} - x_{\min}$ and $J\Delta y = y_{\max} - y_{\min}$, $k \in \{0, \cdots, K\}$, with $K$ such that $K\Delta t = T$.

Most of the Finite Differences in the scheme below are monotone, except for the mixed derivative, which is approximated on a 7-point stencil as in Ma and Forsyth (2014). The expression for the approximation of the mixed derivative depends on the sign of $\rho$. Thus, we distinguish below between the cases $\rho \geq 0$ and $\rho < 0$.

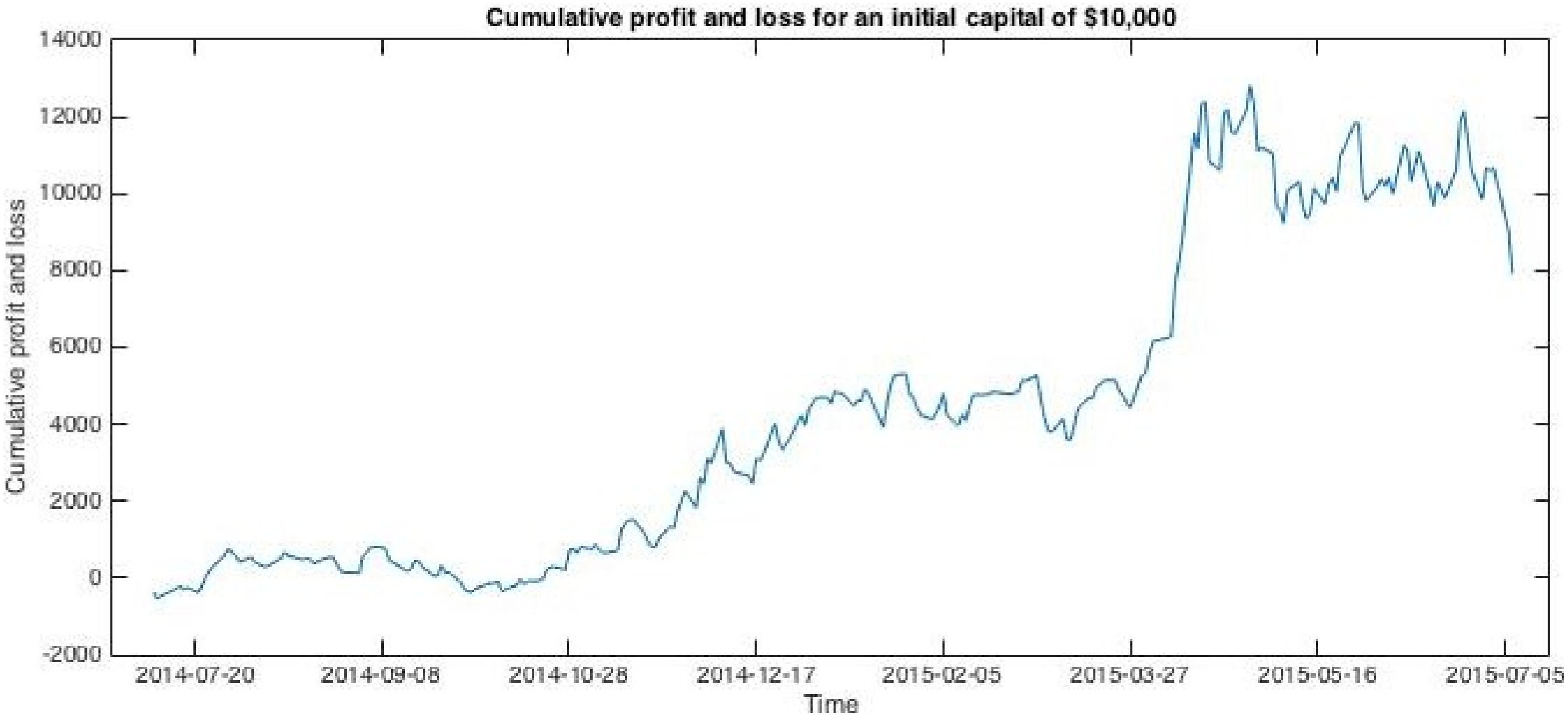

Figure 4.12: This shows the cumulative profit and loss for the out-of-sample test on Bank of China shares, from July 08, 2014 to July 07, 2015. The parameters are estimated from the in-sample interval: $\mu_1 = 0.27$, $\sigma_1 = 0.76$, $\delta_1 = -0.36$, $\theta_1 = -0.16$, $\mu_2 = 0.26$, $\sigma_2 = 0.7$, $\delta_2 = 0.85$, $\theta_2 = -0.37$, $a = -1$, $b = 0.068$, $\beta = -1.12$, $\rho = 0.64$.

For $\rho \geq 0$, the scheme reads

$$
\begin{aligned}
&\frac{\phi_{i,j}^{k+1} - \phi_{i,j}^{k}}{\Delta t} - \frac{\gamma}{2(\gamma-1)^2(1-\rho^2)} \qquad \text{(A.1)}\\
&\left[\frac{\left(\mu_1 - r + \delta_1 z_{i,j}^{k+1}\right)^2}{\sigma_1^2 e^{2\theta_1 x_i}} + \frac{\left(\mu_2 - r + \delta_2 z_{i,j}^{k+1}\right)^2}{\sigma_2^2 e^{2\theta_2 y_j}} - 2\rho \frac{\left(\mu_1 - r + \delta_1 z_{i,j}^{k+1}\right)\left(\mu_2 - r + \delta_2 z_{i,j}^{k+1}\right)}{\sigma_1 \sigma_2 e^{\theta_1 x_i} e^{\theta_2 y_j}}\right] \phi_{i,j}^{k+1}\\
&- \frac{r\gamma}{1-\gamma}\phi_{i,j}^{k+1} + \left(\frac{\mu_1 + \delta_1 z_{i,j}^{k+1} - r\gamma}{\gamma - 1} + \frac{1}{2}\sigma_1^2 e^{2\theta_1 x_i}\right)^+ \frac{\phi_{i,j}^{k+1} - \phi_{i-1,j}^{k+1}}{\Delta x}\\
&- \left(\frac{\mu_1 + \delta_1 z_{i,j}^{k+1} - r\gamma}{\gamma - 1} + \frac{1}{2}\sigma_1^2 e^{2\theta_1 x_i}\right)^- \frac{\phi_{i+1,j}^{k+1} - \phi_{i,j}^{k+1}}{\Delta x}\\
&+ \left(\frac{\mu_2 + \delta_2 z_{i,j}^{k+1} - r\gamma}{\gamma - 1} + \frac{1}{2}\sigma_2^2 e^{2\theta_2 y_j}\right)^+ \frac{\phi_{i,j}^{k+1} - \phi_{i,j-1}^{k+1}}{\Delta y}\\
&- \left(\frac{\mu_2 + \delta_2 z_{i,j}^{k+1} - r\gamma}{\gamma - 1} + \frac{1}{2}\sigma_2^2 e^{2\theta_2 y_j}\right)^- \frac{\phi_{i,j+1}^{k+1} - \phi_{i,j}^{k+1}}{\Delta y}\\
&- \frac{1}{2}\sigma_1^2 e^{2\theta_1 x_i} \frac{\phi_{i-1,j}^{k+1} - 2\phi_{i,j}^{k+1} + \phi_{i+1,j}^{k+1}}{\Delta x^2} - \frac{1}{2}\sigma_2^2 e^{2\theta_2 y_j} \frac{\phi_{i,j-1}^{k+1} - 2\phi_{i,j}^{k+1} + \phi_{i,j+1}^{k+1}}{\Delta y^2}\\
&- \rho\sigma_1\sigma_2 e^{\theta_1 x_i} e^{\theta_2 y_j} \left(\frac{2\phi_{i,j}^{k+1} + \phi_{i+1,j+1}^{k+1} + \phi_{i-1,j-1}^{k+1}}{2\Delta x \Delta y} - \frac{\phi_{i+1,j}^{k+1} + \phi_{i-1,j}^{k+1} + \phi_{i,j+1}^{k+1} + \phi_{i,j-1}^{k+1}}{2\Delta x \Delta y}\right) = 0,
\end{aligned}
$$

and for $\rho < 0$, the scheme reads

$$
\begin{aligned}
&\frac{\phi_{i,j}^{k+1} - \phi_{i,j}^{k}}{\Delta t} - \frac{\gamma}{2(\gamma-1)^2(1-\rho^2)} \qquad \text{(A.2)}\\
&\left[\frac{\left(\mu_1 - r + \delta_1 z_{i,j}^{k+1}\right)^2}{\sigma_1^2 e^{2\theta_1 x_i}} + \frac{\left(\mu_2 - r + \delta_2 z_{i,j}^{k+1}\right)^2}{\sigma_2^2 e^{2\theta_2 y_j}} - 2\rho \frac{\left(\mu_1 - r + \delta_1 z_{i,j}^{k+1}\right)\left(\mu_2 - r + \delta_2 z_{i,j}^{k+1}\right)}{\sigma_1 \sigma_2 e^{\theta_1 x_i} e^{\theta_2 y_j}}\right] \phi_{i,j}^{k+1}
\end{aligned}
$$

$$
\begin{aligned}
&- \frac{r\gamma}{1-\gamma}\phi_{i,j}^{k+1} + \left(\frac{\mu_1 + \delta_1 z_{i,j}^{k+1} - r\gamma}{\gamma - 1} + \frac{1}{2}\sigma_1^2 e^{2\theta_1 x_i}\right)^+ \frac{\phi_{i,j}^{k+1} - \phi_{i-1,j}^{k+1}}{\Delta x} \\
&- \left(\frac{\mu_1 + \delta_1 z_{i,j}^{k+1} - r\gamma}{\gamma - 1} + \frac{1}{2}\sigma_1^2 e^{2\theta_1 x_i}\right)^- \frac{\phi_{i+1,j}^{k+1} - \phi_{i,j}^{k+1}}{\Delta x} \\
&+ \left(\frac{\mu_2 + \delta_2 z_{i,j}^{k+1} - r\gamma}{\gamma - 1} + \frac{1}{2}\sigma_2^2 e^{2\theta_2 y_j}\right)^+ \frac{\phi_{i,j}^{k+1} - \phi_{i,j-1}^{k+1}}{\Delta y} \\
&- \left(\frac{\mu_2 + \delta_2 z_{i,j}^{k+1} - r\gamma}{\gamma - 1} + \frac{1}{2}\sigma_2^2 e^{2\theta_2 y_j}\right)^- \frac{\phi_{i,j+1}^{k+1} - \phi_{i,j}^{k+1}}{\Delta y} nonumber \qquad \text{(A.3)} \\
&- \frac{1}{2}\sigma_1^2 e^{2\theta_1 x_i} \frac{\phi_{i-1,j}^{k+1} - 2\phi_{i,j}^{k+1} + \phi_{i+1,j}^{k+1}}{\Delta x^2} - \frac{1}{2}\sigma_2^2 e^{2\theta_2 y_j} \frac{\phi_{i,j-1}^{k+1} - 2\phi_{i,j}^{k+1} + \phi_{i,j+1}^{k+1}}{\Delta y^2} \\
&- \rho\sigma_1\sigma_2 e^{\theta_1 x_i} e^{\theta_2 y_j} \left(-\frac{2\phi_{i,j}^{k+1} + \phi_{i+1,j-1}^{k+1} + \phi_{i-1,j+1}^{k+1}}{2\Delta x\Delta y} + \frac{\phi_{i+1,j}^{k+1} + \phi_{i-1,j}^{k+1} + \phi_{i,j+1}^{k+1} + \phi_{i,j-1}^{k+1}}{2\Delta x\Delta y}\right) = 0.
\end{aligned}
$$

Furthermore, we impose the Neumann-type boundary conditions

$$
\begin{aligned}
\phi_{1,j}^{k+1} &= \phi_{2,j}^{k+1}, \quad \phi_{I,j}^{k+1} = \phi_{I-1,j}^{k+1}, \\
\phi_{i,1}^{k+1} &= \phi_{i,2}^{k+1}, \quad \phi_{i,J}^{k+1} = \phi_{i,J-1}^{k+1}.
\end{aligned}
$$

Finally the scheme is initialized with

$$
\phi_{i,j}^0 = (1-\gamma)^{\frac{1}{1-\gamma}}.
$$

# Appendix B Generalized Method of Moments estimation

In this section, we describe the econometric approach that we use for estimating the parameters of the stochastic system (2.1), (2.2), and (2.3). We estimate the parameters of the discrete-time discrete dynamics

$$
\ln S_{t+1}^i - \ln S_t^i = \left(\mu_i - \frac{1}{2}\sigma_i^2 e^{2\theta_i \ln S_t^i} + \delta_i \left(a + bt + \ln S_t^1 + \beta \ln S_t^2\right)\right)\Delta t + \epsilon_{t+1}^i, \tag{B.1}
$$

where

$$
\begin{aligned}
\mathbb{E}\left(\epsilon_{t+1}^i\right) &= 0, \\
\mathbb{E}\left[\left(\epsilon_{t+1}^i\right)^2\right] &= \sigma_i^2 e^{2\theta_i \ln S_t^i}\Delta t, \text{ for } i = 1, 2, \\
\mathbb{E}\left[\epsilon_{t+1}^1 \epsilon_{t+1}^2\right] &= \rho\sigma_1 e^{\theta_1 \ln S_t^1}\sigma_2 e^{\theta_2 \ln S_t^2}\Delta t.
\end{aligned}
$$

Please note that, in contrast to the continuous-time model, we no longer require the residuals $\epsilon_{t+1}^i / \left(\sigma_i e^{\theta_i \ln S_t^i}\right)$ to be normally distributed.

We define $\lambda$ to be the parameter vector with elements $\mu_1$, $\sigma_1$, $\delta_1$, $\theta_1$, $\mu_2$, $\sigma_2$, $\delta_2$, $\theta_2$, $\rho$, $a$, $b$, and $\beta$. Given

$$
\epsilon_{t+1}^i = \ln S_{t+1}^i - \ln S_t^i - \left(\mu_i - \frac{1}{2}\sigma_i^2 e^{2\theta_i \ln S_t^i} + \delta_i \left(a + bt + \ln S_t^1 + \beta \ln S_t^2\right)\right)\Delta t, \tag{B.2}
$$

we define the vector $f_t(\lambda)$ of moment functions

$$f_t(\lambda) = \begin{bmatrix} \epsilon^1_{t+1} \\ \epsilon^2_{t+1} \\ \left(\epsilon^1_{t+1}\right)^2 - \sigma_1^2 e^{2\theta_1 \ln S^1_t} \Delta t \\ \left(\epsilon^2_{t+1}\right)^2 - \sigma_2^2 e^{2\theta_2 \ln S^2_t} \Delta t \\ \epsilon^1_{t+1}\epsilon^2_{t+1} - \sigma_1 e^{\theta_1 \ln S^1_t} \sigma_2 e^{\theta_2 \ln S^2_t} \rho \Delta t \end{bmatrix} \bigotimes \begin{bmatrix} 1 \\ t \\ \ln S^1_t \\ \ln S^2_t \end{bmatrix}, \tag{B.3}$$

where $\bigotimes$ stands for the Kronecker Product. Note that $\left(\epsilon^i_{t+1}\right)_i$ plays the role of the unobservable vector of disturbances, while $t$ and $S^i_t$, for $i = 1, 2$, are the chosen instrumental variables. Under the null hypothesis that the restrictions implied by (B.1) are true, we h ve that $\mathbb{E}\left[f_t(\lambda)\right] = 0$.

We then replace $\mathbb{E}\left[f_t(\lambda)\right]$ by its sample counterpart $g_T(\lambda)$

$$g_T(\lambda) = \frac{1}{T} \sum_{t=1}^{T} f_t(\lambda), \tag{B.4}$$

where $T$ is the number of observations and define

$$J_T(\lambda) = g_T^T(\lambda) W_T(\lambda) g_T(\lambda), \tag{B.5}$$

where $g_T(\lambda)$ is given in (B.4), and $W_T(\lambda)$ is a positive-definite symmetric weighting matrix.

Selecting twelve out of the twenty moment functions in (B.3) leads to a model which is just-identified and, in this case, we can simply set the weighting matrix to be the identity matrix. In practice, for the tests reported in this paper, we choose the first eight, followed by the tenth, twelfth, fourteenth, and sixteenth elements.

Beside, in order to assess how well the model fits, we also consider the version of the GMM that uses all of the moment functions in (B.3); the model becomes over-identified, and in order to obtain asymptotically efficient estimates, we set in this case, following Hansen (1982),

$$W_T(\lambda) = S_T^{-1}(\lambda), \tag{B.6}$$

where $S_T^{-1}(\lambda)$ is the estimate of the spectral density matrix proposed by Newey and West (1987), that is,

$$S_T(\lambda) = S_0 + \sum_{j=1}^{k} \left(1 - \frac{j}{k-1}\right) \left(S_j + S_j^T\right), \tag{B.7}$$

where

$$S_j = \frac{1}{T} \sum_{t=j+1}^{T} f_t(\lambda) f_{t-j}^T(\lambda). \tag{B.8}$$

When applying the over-identified model, we estimate the parameters with the following two-step procedure. We first minimize (B.5) by using the identity weighting matrix. Then, we substitute the computed parameters into (B.7) and (B.8) and invert the result to get an approximation of the optimal weighting matrix (B.6). Finally, we minimize (B.5) corresponding to the above approximation of (B.6), in order to compute asymptotically efficient estimates for the parameters.